\documentclass[review]{elsarticle}

\makeatletter
\def\ps@pprintTitle{%
 \let\@oddhead\@empty
 \let\@evenhead\@empty
 \def\@oddfoot{\centerline{\thepage}}%
 \let\@evenfoot\@oddfoot}
\makeatother

\pdfoutput=1
\usepackage{graphicx}
\usepackage{algorithm}
\usepackage{array}
\usepackage[noend]{algpseudocode}
\usepackage{amsfonts, amsmath, amssymb}
\usepackage{stmaryrd}
\usepackage{booktabs,siunitx}
\usepackage[svgnames,table]{xcolor}
\usepackage[tableposition=above]{caption}
\usepackage{pifont}
\usepackage{bm}
\usepackage{cancel}
\usepackage{caption}
\usepackage{color}
\usepackage{enumerate}
\usepackage{float}
\usepackage{hyperref}
\usepackage[english]{babel}
\usepackage[margin=1in]{geometry}
\usepackage{mathtools}
\usepackage{setspace}
\usepackage{subfigure}
\usepackage{cancel}
\usepackage{lineno}

\newcommand \D [2]{\frac{\partial #1}{\partial #2}}

\renewcommand{\vec}[1]{\bm{#1}}
\newcommand{\V}[1]{\bm{#1}}

\def \div{\nabla \cdot \mbox{}}
\def \grad{\nabla}

\def \x{\vec{x}}

\def \u{\vec{u}}

\def \F{\vec{F}}

\def \U{\vec{U}}

\def \F{\vec{F}}

\def \cP{{\mathcal{P}}}

\def \U{\vec{U}}

\def \X{\vec{X}}

\def \cG{\mathcal{G}}

\def \vcL{\vec{\mathcal{L}}}

\def \f{\vec{f}}

\def \3half{\frac{3}{2}}
\def \5half{\frac{5}{2}}

\def \u{\vec{u}}

\def \x{\vec{x}}

\def \div{\nabla \cdot \mbox{}}
\def \grad{\nabla}

\def \cP{{\mathcal{P}}}
\def \InterpSumi{{ \sum_{i = 1}^N}}
\def \InterpSumj{{ \sum_{j = 1}^m}}
\def \cA{\vec{\mathcal{A}}}
\def \cG{\vec{\mathcal{G}}}
\def \cW{\vec{\mathcal{W}}}

\newcommand{\upperRomannumeral}[1]{\uppercase\expandafter{\romannumeral#1}}

\newcommand{\REVIEW}[1]{#1}

\begin{document}

\begin{frontmatter}
	
\title{A unified constraint formulation of  immersed 
body techniques for coupled fluid-solid motion\REVIEW{: continuous equations and numerical algorithms}}
\author[SDSU]{Amneet Pal Singh Bhalla\corref{mycorrespondingauthor}}
\ead{asbhalla@sdsu.edu}
\author[Northwestern1]{Neelesh A. Patankar\corref{mycorrespondingauthor}}
\ead{n-patankar@northwestern.edu}

\address[SDSU]{Department of Mechanical Engineering, San Diego State University, San Diego, CA}
\address[Northwestern1]{Department of Mechanical Engineering, Northwestern University, Evanston, IL}

\cortext[mycorrespondingauthor]{Corresponding author}

\begin{abstract}
Numerical  simulation  of  moving  immersed solid bodies in fluids  is  now  practiced  routinely  following pioneering work of Peskin and co-workers on immersed boundary method (IBM), Glowinski and co-workers  on fictitious  domain  method  (FDM), and others on related methods.  A  variety  of  variants  of  IBM and FDM  approaches have been published, most of which rely on using a background mesh for the fluid equations and tracking the solid body using  Lagrangian  points.  The  key  idea that  is  common  to  these  methods is to  assume  that  the  entire  fluid-solid  domain  is  a  fluid and  then  to  constrain the fluid within the solid domain  to  move in  accordance  with  the solid governing equations. The immersed solid body can be rigid or deforming. Thus, in all these methods the fluid domain is extended into the solid domain\footnote{\REVIEW{Some versions of IBM exclude (or zero-out) the grid points inside the solid domain when solving governing equations. This discussion excludes such methods as they do not fit into the fictitious or extended domain category.}}. In  this  review, we provide a mathemarical perspective of various immersed methods by recasting the governing equations in an extended domain form for the fluid. The solid equations are used to impose appropriate constraints on the fluid that is extended into the solid domain. This leads to extended domain constrained  fluid-solid governing equations that  provide  a  unified  framework  for  various  immersed body techniques. The unified constrained governing equations in the strong form are  independent  of  the  temporal  or  spatial  discretization schemes.  We show that particular choices of  time  stepping  and  spatial  discretization  lead to different techniques reported in literature ranging from freely moving rigid to elastic self-propelling bodies. 
These techniques have wide ranging applications including aquatic  locomotion, underwater  vehicles,  car aerodynamics, and organ physiology (e.g. cardiac flow, esophageal transport, 
respiratory flows), wave energy convertors, among others. We conclude with comments on outstanding challenges and future directions. 

\end{abstract}


\end{frontmatter}

\section{Introduction}
\REVIEW{

The fictitious domain method (FDM) was proposed by Saul'ev in 1963 \cite{saulev1963FDM} to solve diffusion equations on regular meshes. In the 1970s, similar methods emerged known as ``\emph{domain imbedding methods}" where the actual domain is embedded within a regular domain~\cite{buzbee1971DIM}. The immersed boundary method (IBM) is a widely used technique for modeling fluid structure interaction (FSI) following the early work by Peskin \cite{Peskin72b}. IBMs became known in the 1970s and found wider use in 2000s, but their basic ingredients were developed at the end of the 1960s and early 1970s, according to the most recent review article on IBMs by Verzicco~\cite{verzicco2023immersed}. In parallel, the FDM was applied to partial differential equations \cite{astra1978FDM} and in the 1990s, Glowinski and co-workers proposed a distributed Lagrange multiplier (DLM) fictitious domain method for Dirichlet boundary conditions. Later circa 2000 Glowinski and co-workers extended the DLM method to rigid particulate flows \cite{glowinski1999distributed, Patankar2000}.  DLM found broader application in the 2000s by extension to swimming \cite{Shirgaonkar09} and to Brownian systems \cite{sharma2004direct, Chen2006Brownian}. In this review hereafter, all the aforementioned methods and their variants, which rely of extending the fluid domain into the solid domain, will be collectively termed immersed body (IB) methods.

In the past two decades, the IB field has grown tremendously, resulting in several excellent reviews on the topic. In Peskin's 2002~\cite{Peskin02} and Griffith and Patankar's 2020~\cite{griffith2020immersed} papers, IBMs are discussed for modeling biological FSI problems (e.g., cardiovascular and esophageal flows). Mittal and Iaccarino discussed IBMs for flows around rigid bodies in their 2005 review paper~\cite{mittal2005immersed}. Verzicco's 2023 review paper~\cite{verzicco2023immersed} focused on wall modeling approaches to solve turbulent flows. The former two review papers discussed semi-continuous IBM versions, whereas the latter two described fully discrete IBM versions (e.g., direct forcing or velocity reconstruction approaches). These papers provide excellent introductions to IBMs and general guidelines about which IB method is most appropriate for a specific application. However, there is a need to bring different IB methods under the umbrella of a common set of governing equations because of all of them are essentially extended fluid domain implementations. This review fulfills this gap by providing a unified mathematical perspective for various immersed methods. 

Prior presentations and reviews on the immersed boundary method have focused on immersed implementation via the numerical approximations of the continuous equations of motion. In contrast, in this review we focus on the strong form of the governing equations of the extended domain formulation. By numerical approximations, we refer to operations such as interpolation and spreading through delta functions, time-stepping schemes, velocity reconstruction approaches, and operator splitting techniques that work only in time-dependent cases. The strong form of equations, by definition, applies pointwise in the spatial domain, rather than in an average/integral sense. The latter form introduces a numerical length scale into continuous equations, making them semi-continuous instead of fully continuous. The strong form is also valid for steady-state flows (zero Reynolds number) or unsteady (moderate to high Reynolds number) flows and do not rely on any specific numerical treatment. Using particular choices of time stepping and spatial discretization, the strong form of equations can be approximated to produce different IBMs. Consequently, this work contributes to the understanding of the IBM from a theoretical standpoint. This unified mathematical perspective shows that various IBMs are governed by a common, well-defined set of equations that can be analyzed to understand and compare different IB implementations.          

This paper is organized as follows. In Section~\ref{sec_unified_formulation}, we derive the strong form of equations for ``real" fluids and solids, as well as for fictitious (``fake") fluids within real solids. These equations are collectively referred to as fictitious/extended domain equations. In Section~\ref{sec_solution_methods} of the paper, we explain how different choices of time stepping and spatial discretization lead to different IB algorithms. These include the original IBM of Peskin, the immersed finite element method, the velocity forcing or the direct forcing method, the fictitious domain method (FDM) or the distributed Lagrange multiplier method (DLM) for rigid and self-propelling bodies, the Brinkman/volume penalization method, the immersed interface method, and the fully implicit DLM methods. Section~\ref{sec_multiphase_system} provides an overview of some FDM algorithms designed to model three-phase (gas-liquid-solid) FSI problems. Finally, in Section~\ref{sec_miscellaneous} we discuss various issues and improvements to the IB method reported in the literature. Among these are sharp implementations of Lagrange multipliers, stabilization techniques for high-density ratio flows, evaluation of smooth hydrodynamic and constraint forces, incorporating Neumann and Robin boundary conditions on immersed surfaces, fluid leakage prevention via one-sided IB kernels, and overcoming numerical issues arising from narrow gaps in dense particulate flows.       
}

\section{A unified constraint formulation} \label{sec_unified_formulation}

\begin{figure}[]
  \centering
   \includegraphics[scale = 0.4]{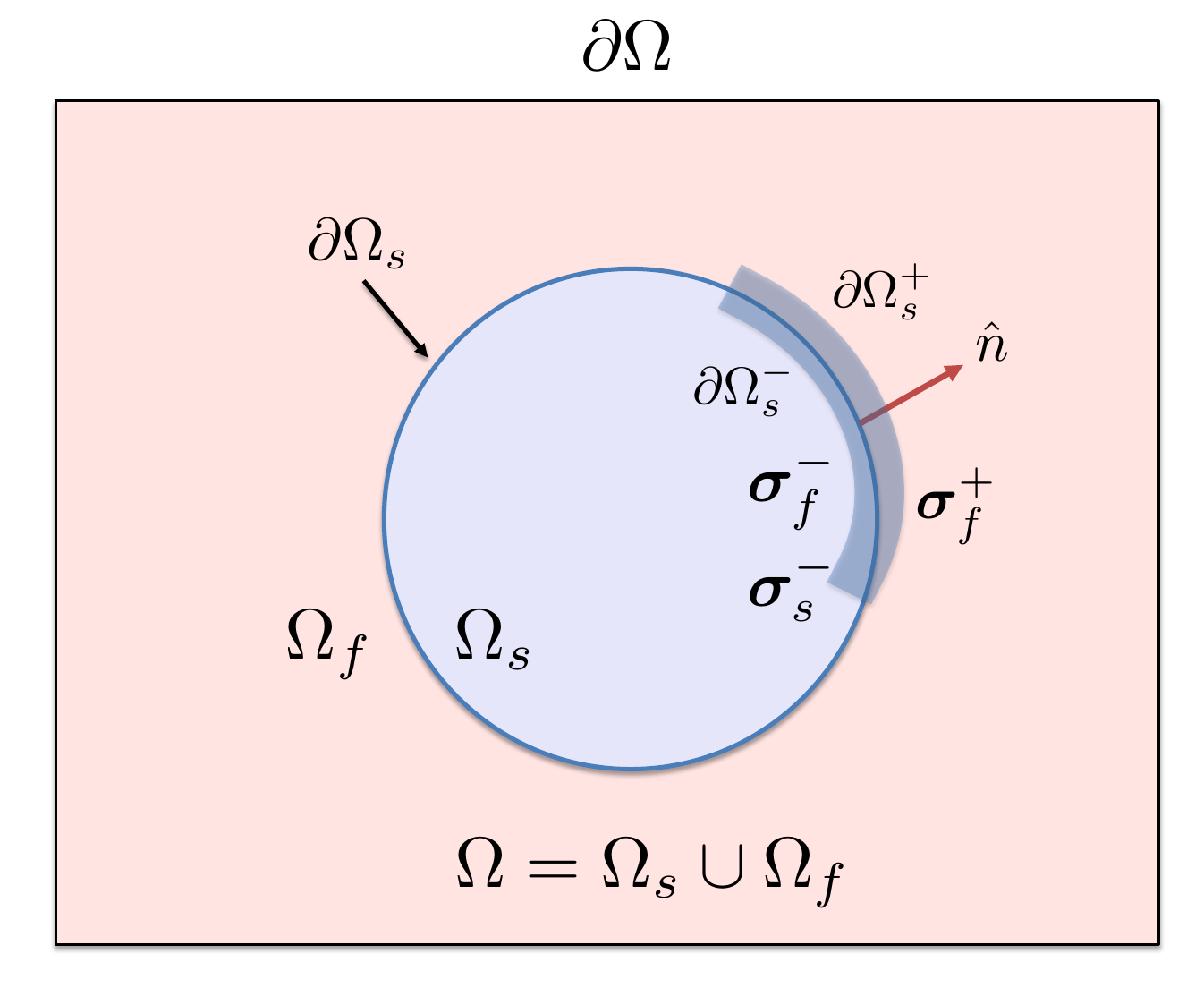}
   \caption{A schematic representation of the fluid-structure interaction system. The computational domain boundary $\partial \Omega$ is demarcated by a solid black line. The filled pink region denotes the fluid domain $\Omega_f$, whereas the filled blue region denotes the solid domain  $\Omega_s$. The solid boundary $\partial \Omega_s$ is shown by a solid blue line. A representative control volume/pillbox used to derive the jump condition given by Equation~\eqref{eq:gov-eqns5} is shown by light grey color. The unit normal vector $\hat{\bm{n}}$ of the solid surface points outwards into the fluid and is shown by a red arrow.}
     \label{fig_FSI_schematic}
\end{figure}

Consider a fluid domain $\Omega_f(t)$ with solid bodies completely immersed in the fluid, as shown in Fig.~\ref{fig_FSI_schematic}. Let the domain of the solid bodies be collectively denoted by  $\Omega_s(t)$. Let the fluid and the solid regions be treated as continua governed by conservation of mass and momentum equations as follows:
\begin{eqnarray}
\label{eq:gov-eqns1}
&&{}\rho_f{D {\bm u}_f\over D t} = {\bm\nabla}\cdot{\bm\sigma}_f~\textrm{in}~\Omega_f(t),\\
\label{eq:gov-eqns2}
&&{}\rho_s{D {\bm u}_s\over D t} = {\bm\nabla}\cdot{\bm\sigma}_s +  {\bm f}_b + \Delta \rho  {\bm g}~\textrm{in}~\Omega_s(t),\\
&&{}{\bm\nabla}\cdot{\bm u}_f=0~\textrm{in}~\Omega_f(t), {\bm\nabla}\cdot{\bm u}_s=0~\textrm{in}~\Omega_s(t),\\
\label{eq:gov-eqns4}
&&{}{\bm u}_f = {\bm u}_s~\textrm{on}~\partial\Omega_s(t),\\
\label{eq:gov-eqns5}
&&{}\biggl[{\bm\sigma}_f^+ - {\bm \sigma}_s^-\biggl]\cdot\hat{\bm n} =- {\bm F}_s~\textrm{on}~\partial\Omega_s(t),
\end{eqnarray}
where $\rho_f$ and $\rho_s$ are the fluid and solid densities, respectively,  ${\bm u}_f$ and ${\bm u}_s$ are the fluid and solid velocity fields, respectively, ${\bm\sigma}_f$ and ${\bm\sigma}_s$ are the fluid and solid stresses, respectively, ${\bm f}_b$ is a body force in the solid other than gravity, $\Delta\rho = \rho_s - \rho_f$, and ${\bm g}$ is the gravitational acceleration. ${D {( )} \over D t}$ denotes a material derivative. Superscript $+$ denotes the value at the interface on fluid side, while superscript $-$ denotes the value at the interface on the solid side. The fluid and the solid are assumed to be incompressible. Note that pressure is included in the stress terms. For simplicity of exposition the hydrostatic pressure in the fluid, due to gravity, is subtracted and no other body force is considered in the fluid. Equations~\ref{eq:gov-eqns4} and~\ref{eq:gov-eqns5} denote the no--slip condition on the fluid-solid interface and the jump in stress, respectively. ${\bm F}_s$ is an external surface force density (i.e. force per unit area) and $\hat{\bm n}$ is an outward normal on the solid surface. In addition, there is a boundary condition \REVIEW{of Dirichlet type}\footnote{\REVIEW{Traction or natural boundary conditions can also be considered on (part of) the domain boundary. This will lead to some additional terms in the weak formulation. The strong form of the equations, however, will remain the same. Dirichlet boundary conditions are considered here to simplify the strong form derivation.}} for velocity at the boundary of the entire domain $\Omega$, and there are initial conditions for the fluid and solid velocities, which are not listed above. 

Equations~\ref{eq:gov-eqns1}--\ref{eq:gov-eqns5} are the strong form equations for the fluid and solid domains. The corresponding combined weak form of the governing equations, where the interface conditions have been used, are given below
\begin{eqnarray}
\label{eq:gov-eqns-weak}
&&{}\int\limits_{\Omega_f}\rho_f{D {\bm u}_f\over D t}\cdot{\bm v}_f dV
 +\int\limits_{\Omega_f}{\bm\sigma}_f:{\bm D}({\bm v}_f)dV 
 -\int\limits_{\Omega_f}q_f{\bm\nabla}\cdot{\bm u}_fdV \nonumber \\
 &&{}+\int\limits_{\Omega_s}\rho_s{D {\bm u}_s\over D t}\cdot{\bm v}_s dV
 +\int\limits_{\Omega_s}{\bm\sigma}_s:{\bm D}({\bm v}_s)dV 
 -\int\limits_{\Omega_s}q_s{\bm\nabla}\cdot{\bm u}_sdV \nonumber \\
&&{}-\int\limits_{\Omega_s}({\bm f}_b + \Delta \rho {\bm g} ) \cdot {\bm v}_s dV
-\int\limits_{\partial\Omega_s}{\bm v}_s \cdot {\bm F}_s dS \nonumber \\
&&{}=0\;\;\forall\; {\bm v}_f\in\;S_{fv},\forall\; {\bm v}_s\in\;S_{sv},\;\forall \; q_f\in L^2(\Omega_f(t)),\;\forall \; q_s\in L^2(\Omega_s(t)),
\end{eqnarray}
where ${\bm v}_f$ is the variation of ${\bm u}_f$, ${\bm v}_s$ is the variation of ${\bm u}_s$, and $q_f/q_s$ are the variations of pressure $p_f/p_s$ (which are within the stress terms). The operator ${\bm D}()$ acts on a vector, say ${\bm u}$, as follows: ${\bm D}(\bm u) = {1 \over 2}({\bm\nabla}{\bm u} + {\bm\nabla}{\bm u}^T)$. The solution spaces of pressure are the same as their variances. The solution spaces $S_f$ (for ${\bm u}_f$) and $S_s$ (for ${\bm u}_s$), and the variation spaces $S_{fv}$ (for ${\bm v}_f$) and $S_{sv}$ (for ${\bm v}_s$) are given by
\begin{eqnarray}
\label{eq:spaces1}
&&{}S_f = \{{\bm u}_f | {\bm u}_f \in H^{1}(\Omega_f(t))^{3}, {\bm u}_f = {\bm u}_{\partial\Omega}~\textrm{on}~\partial\Omega, {\bm u}_f = {\bm u}_s~\textrm{on}~\partial\Omega_s(t) \},\\
&&{}S_s = \{{\bm u}_s | {\bm u}_s \in H^{1}(\Omega_s(t))^{3}, {\bm u}_s = {\bm u}_f~\textrm{on}~\partial\Omega_s(t) \},\\
&&{}S_{fv} = \{{\bm v}_f | {\bm v}_f \in H^{1}(\Omega_f(t))^{3}, {\bm v}_f = 0~\textrm{on}~\partial\Omega, {\bm v}_f = {\bm v}_s~\textrm{on}~\partial\Omega_s(t) \},\\
&&{}S_{sv} = \{{\bm v}_s | {\bm v}_s \in H^{1}(\Omega_s(t))^{3}, {\bm v}_s = {\bm v}_f~\textrm{on}~\partial\Omega_s(t) \}.
\end{eqnarray}

The key feature of all immersed techniques is to extend the fluid equations into the solid domain.  The motivation to do so is the convenience of solving the fluid equations in the entire fluid-solid domain without having to worry about the continuously changing fluid domain and keeping track of it. The additional (but redundant) solution of the fluid equations in the solid domain should not affect the original solution of fluid-solid motion. Thus, we note that the fluid equations extended into the solid domain should be solved such that the velocity on the fluid-solid interface $\partial\Omega_s(t)$ is ``known" from the solution of fluid-solid motion. This can be done in weak form equations by extending the fluid equations into the solid domain as follows
\begin{eqnarray}
\label{eq:extd-dom-weak}
&&{}\int\limits_{\Omega_s}\rho_f{D {\bm u}_f\over D t}\cdot{\bm v}_e dV
 +\int\limits_{\Omega_s}{\bm\sigma}_f:{\bm D}({\bm v}_e)dV 
 -\int\limits_{\Omega_s}q_f{\bm\nabla}\cdot{\bm u}_f dV \nonumber \\
&&=0\;\;\forall\; {\bm v}_e = {\bm v}_f - {\bm v}_s, \;\forall \; q_f\in L^2(\Omega_s(t)),
\end{eqnarray}
where ${\bm v}_e$ is the variation of the fluid velocity field ${\bm u}_f$ extended into the solid domain. We need ${\bm v}_e\in H^{1}(\Omega_s(t))^{3}$ such that ${\bm v}_e=0$ on $\partial\Omega_s(t)$; the latter condition arising because the extended velocity field ${\bm u}_f$ must equal the ``known" velocity at the fluid-solid interface as discussed above. Without loss of generality we choose ${\bm v}_e = {\bm v}_f - {\bm v}_s$, where ${\bm v}_f$ is the variation of the extended fluid velocity field. This choice satisfies all the requirements on ${\bm v}_e$. The spaces for ${\bm u}_f$ and ${\bm v}_f$, after accounting for the extended domain, will be formally identified below. 

Combining Equations \ref{eq:gov-eqns-weak} and \ref{eq:extd-dom-weak}, the following extended domain weak form of the governing equations is obtained
\begin{eqnarray}
\label{eq:combined-weak}
&&{}\int\limits_{\Omega}\rho_f{D {\bm u}_f\over D t}\cdot{\bm v}_f dV
 +\int\limits_{\Omega}{\bm\sigma}_f:{\bm D}({\bm v}_f)dV 
 -\int\limits_{\Omega}q_f{\bm\nabla}\cdot{\bm u}_fdV \nonumber \\
 &&{}+\int\limits_{\Omega_s} \biggl( \rho_s {D {\bm u}_s\over D t} - \rho_f {D {\bm u}_f\over D t} \biggl) \cdot{\bm v}_s dV
 +\int\limits_{\Omega_s}\Delta {\bm\sigma}:{\bm D}({\bm v}_s)dV 
 -\int\limits_{\Omega_s}q_s{\bm\nabla}\cdot{\bm u}_sdV \nonumber \\
&&{}-\int\limits_{\Omega_s}({\bm f}_b + \Delta \rho {\bm g} ) \cdot {\bm v}_s dV
-\int\limits_{\partial\Omega_s}{\bm v}_s \cdot {\bm F}_s dS \nonumber \\ 
&&{}=0\;\;\forall\; {\bm v}_f\in\ S_{ev},\forall\; {\bm v}_s\in\ S_{sv},\;\forall \; q_f\in L^2(\Omega),\;\forall \; q_s\in L^2(\Omega_s(t)),
\end{eqnarray}
where $\Delta{\bm\sigma} = {\bm\sigma}_s - {\bm\sigma}_f$ and ${\bm v}_e = {\bm v}_f - {\bm v}_s$ has been used.  The solution space $S_e$ (for ${\bm u}_f$) and the variation space $S_{ev}$ (for ${\bm v}_f$) are given by
\begin{eqnarray}
\label{eq:ext-spaces1}
&&{}S_e = \{{\bm u}_f | {\bm u}_f \in H^{1}(\Omega)^{3}, {\bm u}_f = {\bm u}_{\partial\Omega}~\textrm{on}~\partial\Omega, {\bm u}_f = {\bm u}_s~\textrm{on}~\partial\Omega_s(t) \},\\
&&{}S_{ev} = \{{\bm v}_f | {\bm v}_f \in H^{1}(\Omega)^{3}, {\bm v}_f = 0~\textrm{on}~\partial\Omega, {\bm v}_f = {\bm v}_s~\textrm{on}~\partial\Omega_s(t) \}.
\end{eqnarray}

The next step is to relax the constraints (${\bm u}_f = {\bm u}_s$ and ${\bm v}_f = {\bm v}_s$ on $\partial\Omega_s(t)$) in solution and variation spaces of the velocity fields. This is compensated by enforcing the constraints on velocities and adding the corresponding Lagrange multipliers in the governing equations. This can be done in different ways. Two specific formulations, which will be called the body force and stress formulations are presented below.

\subsection{The body force formulation}

The constraint that should necessarily be imposed for an extended fluid domain formulation is ${\bm u}_f = {\bm u}_s$ on $\partial\Omega_s(t)$. In this case the extended fluid and solid velocities will necessarily match only on the fluid-solid interface. Another legitimate option is to require that the extended fluid velocity field be equal to the solid velocity field in the entire solid domain (including the fluid-solid interface). We will refer to this latter constraint as imposing ${\bm u}_f = {\bm u}_s$ on $\mathcal{B}(t)$, where $\mathcal{B}(t)$ represents the entire solid domain including the fluid-solid interface. ${\bm u}_f = {\bm u}_s$ on $\mathcal{B}(t)$ can be imposed in the weak form as follows
\begin{eqnarray}
\label{eq:constraint-body}
&&{} \int\limits_{\partial\Omega_s} ( {\bm u}_s - {\bm u}_f) \cdot  {\bm \Psi_s} dS +\int\limits_{\Omega_s} ( {\bm u}_s - {\bm u}_f) \cdot  {\bm \psi_b} dV = 0~\forall\; {\bm \Psi_s}\in L^2(\partial\Omega_s(t))^3, \forall\; {\bm \psi_b}\in L^2(\Omega_s(t))^3.
\end{eqnarray}
If one chooses to impose ${\bm u}_f = {\bm u}_s$ on $\partial\Omega_s(t)$, then only the first integral in the equation above is set to zero; the second integral is not present. In summary, the first integral is essential while the second integral may or may not be used. 

It is also necessary to add an appropriate distributed Lagrange multiplier term, corresponding to the above constraint, in the weak form of the momentum equation. The extended domain weak form (Equation \ref{eq:combined-weak}), together with the constraint (Equation \ref{eq:constraint-body}) and the corresponding Lagrange multiplier terms, becomes
\begin{eqnarray}
\label{eq:combined-extd-dom-weak}
&&{}\int\limits_{\Omega}\rho_f{D {\bm u}_f\over D t}\cdot{\bm v}_f dV
 +\int\limits_{\Omega}{\bm\sigma}_f:{\bm D}({\bm v}_f)dV 
 -\int\limits_{\Omega}q_f{\bm\nabla}\cdot{\bm u}_fdV \nonumber \\
 &&{}+\int\limits_{\Omega_s} \biggl( \rho_s {D {\bm u}_s\over D t} - \rho_f {D {\bm u}_f\over D t} \biggl) \cdot{\bm v}_s dV
 +\int\limits_{\Omega_s} \Delta{\bm\sigma}:{\bm D}({\bm v}_s)dV 
 -\int\limits_{\Omega_s}q_s{\bm\nabla}\cdot{\bm u}_sdV \nonumber \\
&&{}+\int\limits_{\partial\Omega_s} ( {\bm u}_s - {\bm u}_f) \cdot  {\bm \Psi_s} dS 
+\int\limits_{\Omega_s} ( {\bm u}_s - {\bm u}_f) \cdot  {\bm \psi_b} dV
+\int\limits_{\partial\Omega_s} ( {\bm v}_s - {\bm v}_f) \cdot  {\bm \Lambda_s} dS  
+\int\limits_{\Omega_s} ( {\bm v}_s - {\bm v}_f) \cdot  {\bm \lambda_b} dV \nonumber \\
&&{}-\int\limits_{\Omega_s}({\bm f}_b + \Delta \rho {\bm g} ) \cdot {\bm v}_s dV
-\int\limits_{\partial\Omega_s}{\bm v}_s \cdot {\bm F}_s dS = 0~\forall\; {\bm v}_f\in S_{ev0},\forall\; {\bm v}_s\in H^1(\Omega_s(t))^3, \nonumber \\ 
&&{}\;\forall \; q_f\in L^2(\Omega), \forall \; q_s\in L^2(\Omega_s(t)), \forall\; {\bm \Psi_s}\in L^2(\partial\Omega_s(t))^3, \forall\; {\bm \psi_b}\in L^2(\Omega_s(t))^3,
\end{eqnarray}
where ${\bm u}_s\in H^1(\Omega_s(t))^3$, $\bm \Lambda_s\in L^2(\partial\Omega_s(t))^3$ is the distributed Lagrange multiplier (DLM)  field corresponding to the first integral in Equation \ref{eq:constraint-body} (the essential surface constraint) and $\bm \lambda_b\in L^2(\Omega_s(t))^3$ is the (distributed) Lagrange multiplier field corresponding to the second integral in Equation \ref{eq:constraint-body}. It will be seen below that $\bm \Lambda_s$ emerges as a force per unit area on the fluid-solid interface, whereas $\bm \lambda_b$ emerges as a force per unit volume in the solid domain. The solution space $S_{e0}$ (for ${\bm u}_f$) and the variation space $S_{ev0}$ (for ${\bm v}_f$) are given by
\begin{eqnarray}
\label{eq:ext2-spaces1}
&&{}S_{e0} = \{{\bm u}_f | {\bm u}_f \in H^{1}(\Omega)^{3}, {\bm u}_f = {\bm u}_{\partial\Omega}~\textrm{on}~\partial\Omega \},\\
&&{}S_{ev0} = \{{\bm v}_f | {\bm v}_f \in H^{1}(\Omega)^{3}, {\bm v}_f = 0~\textrm{on}~\partial\Omega \}.
\end{eqnarray}
In order to obtain the extended domain strong form, Equation \ref{eq:combined-extd-dom-weak} is converted to the following form after using the Gauss theorem
\begin{eqnarray}
\label{eq:combined-extd-dom-weak-reorg}
&&{}\int\limits_{\Omega} \biggl( \rho_f{D {\bm u}_f\over D t} - {\bm\nabla}\cdot{\bm\sigma}_f - {\bm \Lambda_s} \delta_s - {\bm \lambda_b} \biggl)\cdot {\bm v}_f dV \nonumber \\
 &&{}+\int\limits_{\Omega_s} \biggl( \rho_s{D {\bm u}_s\over D t} - \rho_f{D {\bm u}_f\over D t} 
 - {\bm\nabla}\cdot\Delta{\bm\sigma}-{\bm f_b}- \Delta\rho{\bm g} +{\bm \lambda}_b \biggl) \cdot {\bm v}_sdV \nonumber \\
&&{} -\int\limits_{\Omega}q_f{\bm\nabla}\cdot{\bm u}_fdV  -\int\limits_{\Omega_s}q_s{\bm\nabla}\cdot{\bm u}_sdV +\int\limits_{\partial\Omega_s} ( {\bm u}_s - {\bm u}_f) \cdot  {\bm \Psi_s} dS +\int\limits_{\Omega_s} ( {\bm u}_s - {\bm u}_f) \cdot  {\bm \psi_b} dV \nonumber \\
&&{}+\int\limits_{\partial\Omega_s}{\bm v}_s \cdot (\Delta {\bm \sigma} \cdot \hat {\bm n}+{\bm \Lambda}_s-{\bm F}_s) dS =0~\forall\; {\bm v}_f\in S_{ev0},\forall\; {\bm v}_s\in H^1(\Omega_s(t))^3,  \nonumber \\
&&{}\;\forall \; q_f\in L^2(\Omega), \forall \; q_s\in L^2(\Omega_s(t)), \forall\; {\bm \Psi_s}\in L^2(\partial\Omega_s(t))^3, \forall\; {\bm \psi_b}\in L^2(\Omega_s(t))^3,
\end{eqnarray}
where $\delta_s$ is a surface delta function, i.e., non-zero only on $\partial\Omega_s(t)$. This leads to the following extended domain strong form
\begin{eqnarray}
\label{eq:fluid-mom}
&&{}\rho_f{D {\bm u}_f\over D t} = {\bm\nabla}\cdot{\bm\sigma}_f+{\bm \Lambda_s} \delta_s + {\bm \lambda_b}~\textrm{in}~\Omega, \\
\label{eq:solid-mom}
&&{}\rho_s{D {\bm u}_s\over D t} - \rho_f{D {\bm u}_f\over D t} = {\bm\nabla}\cdot\Delta{\bm\sigma}+{\bm f_b} + \Delta\rho{\bm g}-{\bm \lambda_b}~\textrm{in}~\Omega_s(t), \\
\label{eq:cont}
&&{}{\bm\nabla}\cdot{\bm u}_f=0~\textrm{in}~\Omega, \quad {\bm\nabla}\cdot{\bm u}_s=0~\textrm{in}~\Omega_s(t),\\
\label{eq:vel-constraint}
&&{}{\bm u}_f = {\bm u}_s~\textrm{on}~\mathcal{B}(t), \\
\label{eq:solid-mom-BC}
&&{}\Delta\bm\sigma\cdot\hat{\bm n} =
\biggl[{\bm\sigma}_s^- - {\bm \sigma}_f^-\biggl]\cdot\hat{\bm n} = {\bm F}_s - {\bm \Lambda}_s~\textrm{on}~\partial\Omega_s(t).
\end{eqnarray}

Equation \ref{eq:fluid-mom} is the momentum equation of the fluid extended to the entire domain. The forces $\bm \Lambda_s$ and $\bm \lambda_b$ in the fluid equation are non-zero only on the fluid-solid interface and the solid domain, respectively. These forces arise due to the constraint that the extended fluid in the solid domain moves with the solid body (Equation \ref{eq:vel-constraint}). Specifically, $\bm \Lambda_s$ (force per unit area) arises because of the constraint ${\bm u}_f = {\bm u}_s~\textrm{on}~\partial\Omega_s(t)$ (which is included in ${\bm u}_f = {\bm u}_s~\textrm{on}~\mathcal{B}(t)$) and it should always be present in a formally correct formulation. $\bm \lambda_b$ (force per unit volume) arises only if the constraint ${\bm u}_f = {\bm u}_s$ is imposed on the entire solid domain; if ${\bm u}_f = {\bm u}_s$ is imposed only on $\partial\Omega_s(t)$ then ${\bm \lambda_b} = {\bm 0}$. Note that the fluid velocity field outside the body and the solid velocity field will be the same correct solution irrespective of whether ${\bm u}_f = {\bm u}_s~\textrm{on}~\partial\Omega_s(t)$ is imposed or ${\bm u}_f = {\bm u}_s~\textrm{on}~\mathcal{B}(t)$ is imposed; only the extended fluid velocity field inside the solid domain differs in these scenarios (which is not a solution of physical interest anyway).

Equation \ref{eq:solid-mom} is the correction to the momentum equation in the solid domain. In other words, adding Equations \ref{eq:fluid-mom} and \ref{eq:solid-mom} in $\Omega_s(t)$ gives the momentum equation of the solid body (Equation \ref{eq:gov-eqns2}). Equation \ref{eq:solid-mom-BC} acts as a boundary condition for  Equation \ref{eq:solid-mom}. The implication of this boundary condition becomes clear as follows. Consider a pillbox (control volume of zero thickness) at any location on the fluid-solid interface and apply Equation \ref{eq:fluid-mom} to the pillbox to get
\begin{equation}
\label{eq:fluid-jump}
\biggl[{\bm\sigma}_f^+ - {\bm \sigma}_f^-\biggl]\cdot\hat{\bm n} = -{\bm \Lambda}_s~\textrm{on}~\partial\Omega_s(t).
\end{equation}
This implies that the surface force density $\bm \Lambda_s \delta_s$ causes a jump in fluid stress across the fluid-solid interface. Subtracting Equation \ref{eq:solid-mom-BC} from Equation \ref{eq:fluid-jump} leads to the stress jump condition (Equation \ref{eq:gov-eqns5}) in the original strong form of this problem. Finally, Equation \ref{eq:cont} imposes the incompressibility constraint that leads to the fluid and solid pressures (which are within the stress terms). 

Equations \ref{eq:fluid-mom}--\ref{eq:solid-mom-BC} represent the body force form of the extended domain constraint formulation of fluid-solid equations. In the next section we consider the stress formulation of the same problem. Equations \ref{eq:fluid-mom}--\ref{eq:solid-mom-BC} are different but equivalent to the original strong form in Equations \ref{eq:gov-eqns1}--\ref{eq:gov-eqns5}. The difference arises due to the extension of the fluid into the solid domain, which is not the case in the original form. Equations \ref{eq:fluid-mom}--\ref{eq:solid-mom-BC} can be solved by different algorithms, which have been known in literature as immersed boundary methods (IBM) or fictitious domain methods (FDM), among others. This will be elaborated in a later section. 

\subsection{The stress formulation}

The constraint in Equation \ref{eq:constraint-body} can be imposed by using different inner products. This can lead to different formulations. The following choice of inner product leads to what we call the stress formulation
\begin{eqnarray}
\label{eq:constraint-stress}
&&{} \int\limits_{\partial\Omega_s} ( {\bm u}_s - {\bm u}_f) \cdot  {\bm \Psi_s} dS -\int\limits_{\Omega_s} {\bm D}({\bm u}_s - {\bm u}_f):{\bm D}({\bm\xi_b}) dV \nonumber \\
&&{}= 0~\forall\; {\bm \Psi_s}\in L^2(\partial\Omega_s(t))^3, \forall\; {\bm \xi_b}\in H^1(\Omega_s(t))^3.
\end{eqnarray}
Equation \ref{eq:constraint-stress} imposes ${\bm u}_f - {\bm u}_s$ on $\mathcal{B}(t)$. This will become evident after the strong form is derived below. As discussed in the previous section, if only ${\bm u}_f - {\bm u}_s$ on $\partial\Omega_s(t)$ is to be imposed, then the second integral in the equation above will not be present. The difference between body force and stress formulations arise only when ${\bm u}_f - {\bm u}_s$ on $\mathcal{B}(t)$ is imposed. If ${\bm u}_f - {\bm u}_s$ on $\partial\Omega_s(t)$ is imposed, then both formulations are identical. 

The extended domain weak form (Equation \ref{eq:combined-weak}), together with the constraint (Equation \ref{eq:constraint-stress}) and the corresponding Lagrange multiplier terms, becomes
\begin{eqnarray}
\label{eq:combined-extd-dom-stress-weak}
&&{}\int\limits_{\Omega}\rho_f{D {\bm u}_f\over D t}\cdot{\bm v}_f dV
 +\int\limits_{\Omega}{\bm\sigma}_f:{\bm D}({\bm v}_f)dV -\int\limits_{\Omega}q_f{\bm\nabla}\cdot{\bm u}_fdV \nonumber \\
&&{}+\int\limits_{\Omega_s} \biggl( \rho_s {D {\bm u}_s\over D t} - \rho_f {D {\bm u}_f\over D t} \biggl) \cdot{\bm v}_s dV
 +\int\limits_{\Omega_s} \Delta{\bm\sigma}:{\bm D}({\bm v}_s)dV 
-\int\limits_{\Omega_s}q_s{\bm\nabla}\cdot{\bm u}_sdV \nonumber \\
&&{}+\int\limits_{\partial\Omega_s} ({\bm u}_s - {\bm u}_f) \cdot  {\bm \Psi_s} dS -\int\limits_{\Omega_s}{\bm D}({\bm u}_s - {\bm u}_f):{\bm D}({\bm\xi_b})dV +\int\limits_{\partial\Omega_s} ( {\bm v}_s - {\bm v}_f) \cdot  {\bm \Lambda}_s dS -\int\limits_{\Omega_s}{\bm D}({\bm v}_s - {\bm v}_f):{\bm D}({\bm\zeta_b})dV \nonumber \\
&&{}-\int\limits_{\Omega_s}({\bm f}_b + \Delta \rho {\bm g} ) \cdot {\bm v}_s dV -\int\limits_{\partial\Omega_s}{\bm v}_s \cdot {\bm F}_s dS =0~\forall\; {\bm v}_f\in S_{ev0},\forall\; {\bm v}_s\in H^1(\Omega_s(t))^3,  \nonumber \\
&&{}\;\forall \; q_f\in L^2(\Omega), \forall \; q_s\in L^2(\Omega_s(t)), \forall\; {\bm \Psi_s}\in L^2(\partial\Omega_s(t))^3, \forall\; {\bm \xi_b}\in H^1(\Omega_s(t))^3,
\end{eqnarray}
where $\bm \zeta_b\in H^1(\Omega_s(t))^3$ is the Lagrange multiplier field corresponding to the second integral in Equation \ref{eq:constraint-stress}. All other solution and variation spaces are the same as in the body force formulation. Equation \ref{eq:combined-extd-dom-stress-weak} is converted to another weak form after using the Gauss theorem as was done earlier for the body force formulation. The resulting extended domain strong form is 
\begin{eqnarray}
\label{eq:fluid-stress-mom}
&&{}\rho_f{D {\bm u}_f\over D t} = {\bm\nabla}\cdot{\bm\sigma}_f
+{\bm\nabla}\cdot{\bm D}({\bm\zeta_b})+{\bm \Lambda_s}\delta_s~\textrm{in}~\Omega, \\
\label{eq:solid-stress-mom}
&&{}\rho_s{D {\bm u}_s\over D t} - \rho_f{D {\bm u}_f\over D t} = {\bm\nabla}\cdot\Delta{\bm\sigma}+{\bm f_b} + \Delta\rho{\bm g}-{\bm\nabla}\cdot{\bm D}({\bm\zeta_b})~\textrm{in}~\Omega_s(t), \\
\label{eq:vel-stress-body-constraint}
&&{}{\bm\nabla}\cdot({\bm D}({\bm u}_f - {\bm u}_s))={\bm 0}~\textrm{in}~\Omega_s(t);\quad
{\bm D}({\bm u}_f - {\bm u}_s) \cdot \hat{\bm n}={\bm 0}~\textrm{on}~\partial\Omega_s(t),\\
\label{eq:vel-stress-constraint}
&&{}{\bm u}_f - {\bm u}_s = {\bm 0}
\quad\textrm{on}\;\; \partial\Omega_s(t), \\
\label{eq:solid-stress-mom-BC}
&&{}\biggl[{\bm\sigma}_s^- - {\bm \sigma}_f^- - {\bm\sigma}_{\zeta}^- \biggl]\cdot\hat{\bm n} = {\bm F}_s-{\bm \Lambda}_s~\textrm{on}~\partial\Omega_s(t),
\end{eqnarray}
where ${\bm\sigma}_{\zeta} = {\bm D}({\bm\zeta_b})$ and the incompressibility constraint equations are not repeated here (Equation \ref{eq:cont}). Equations \ref{eq:fluid-stress-mom}--\ref{eq:solid-stress-mom-BC} represent the stress form of the extended domain constraint formulation of fluid-solid equations. As in the body force formulation the Lagrange multiplier force field $\bm \Lambda_s$ arises due to the surface constraint Equation \ref{eq:vel-stress-constraint}. However, the other Lagrange multiplier $\bm \zeta_b$ has a different form compared to $\bm \lambda_b$ in the body force formulation. ${\bm\sigma}_{\zeta}= {\bm D}({\bm\zeta_b})$ is the stress field corresponding to the Lagrange multiplier $\bm \zeta_b$. $\bm \zeta_b$ arises because of the constraint in Equation \ref{eq:vel-stress-body-constraint}, which ensures that the extended fluid velocity field in the solid domain can differ from the solid velocity field in that region only by a rigid motion. This rigid motion is constrained to be zero by Equation \ref{eq:vel-stress-constraint}.  Thus, Equations \ref{eq:vel-stress-body-constraint} and \ref{eq:vel-stress-constraint} collectively ensure that ${\bm u}_f - {\bm u}_s$ on $\mathcal{B}(t)$. 

As before, adding Equations \ref{eq:fluid-stress-mom} and \ref{eq:solid-stress-mom} in $\Omega_s(t)$ gives the momentum equation of the solid body (Equation \ref{eq:gov-eqns2}). Equation \ref{eq:solid-stress-mom-BC} acts as a boundary condition for  Equation \ref{eq:solid-stress-mom}. Equation \ref{eq:fluid-stress-mom} implies the following jump condition for fluid stress across the fluid-solid interface
\begin{equation}
\label{eq:fluid-stress-jump}
\biggl[{\bm\sigma}_f^+ - {\bm \sigma}_f^- - {\bm\sigma}_{\zeta}^- \biggl]\cdot\hat{\bm n} = -{\bm \Lambda}_s~\textrm{on}~\partial\Omega_s(t).
\end{equation}
Subtracting Equation \ref{eq:solid-stress-mom-BC} from Equation \ref{eq:fluid-stress-jump} leads to the stress jump condition (Equation \ref{eq:gov-eqns5}) in the original strong form of this problem.

\section{Solution methods and algorithms} \label{sec_solution_methods}
Several families of computational methods, where the fluid domain is extended into the solid domain, have been developed to solve the fluid-solid motion problem. The governing equations being solved in all these methods are fundamentally the same as those presented in the previous section. Different methods are obtained depending on how the governing equations are solved. We will elucidate this below by considering the body force form of the governing equations (Equations \ref{eq:fluid-mom}--\ref{eq:solid-mom-BC}). We will avoid spatial discretization so that the key algorithmic ideas, that define the different methods, are evident. Temporal discretization will also be minimal in our discussion; elementary fractional time-stepping (FTS) schemes will be used in some cases for clarity. The discretization choices are important in implementing these algorithms but that is not the focus of the discussion below. The incompressibility constraint is implied in the equations below and will not be explicitly mentioned.

\subsection{Immersed boundary method} \label{sec_ibm}

Peskin's immersed boundary method \cite{Peskin02} is one of the most widely used techniques to simulate coupled fluid-elastic-body motion. The original form considered elastic fibers immersed in a fluid \cite{Peskin72b}. The approach has since been generalized to simulate 3D elastic continua immersed in a fluid \cite{HGao14-MV,HGao14-iblv_diastole}. When Equations \ref{eq:fluid-mom}--\ref{eq:solid-mom-BC} are solved in a particular order, it is equivalent to the classical IBM.
  \begin{algorithm}[]
  \caption{IMMERSED BOUNDARY METHOD}\label{IBM-algo}
  \begin{algorithmic}
    \State \textbf{1.} Solve  \\ \\
    $\quad \quad \rho_f {D {\bm u}_f\over D t} = {\bm\nabla}\cdot{\bm\sigma}_f+{\bm \lambda}$ in $\Omega$. \Comment{obtain the velocity field $\V{u}_f$ in the entire domain} \\
    
   \State \textbf{2.} Impose the constraint ${\bm u}_s = {\bm u}_f~\textrm{on}~\mathcal{B}(t)$. \Comment{e.g. use regularized kernels to interpolate $\V{u}_f$ onto $\mathcal{B}(t)$} \\

   \State \textbf{3.} Compute \\ \\
   $\quad \quad {\bm \lambda_b} = -\Delta \rho {D {\bm u}_s\over D t}  + {\bm\nabla}\cdot({{\bm\sigma}_s - {\bm\sigma}_f})+{\bm f_b} + \Delta\rho{\bm g}~\textrm{in}~\Omega_s(t)$, \\ \\
   $\quad \quad {\bm \Lambda}_s = -\biggl[{\bm\sigma}_s^- - {\bm \sigma}_f^-\biggl]\cdot\hat{\bm n} + {\bm F}_s ~\textrm{on}~\partial\Omega_s(t)$, \\ \\
   $\quad \quad {\bm \lambda} = {\bm \Lambda}_s \delta_s + {\bm \lambda}_b$. \Comment{update the Lagrange multipliers} \\

\State \textbf{4.} Go to step 1.

  \end{algorithmic}
  \end{algorithm}

Typically the calculation at each time-step of an IBM proceeds in three explicit steps as described in Algorithm~\ref{IBM-algo}. In step 1, Equation \ref{eq:fluid-mom} from the body force formulation is solved with the latest known $\bm \lambda$ field. Note that $\bm \Lambda_s$ is typically smeared over few grids at the fluid-solid interface using discrete delta functions \cite{Peskin02}. Thus, in practice the entire $\bm \lambda$ field is treated like a body force term. In step 2, the constraint in Equation \ref{eq:vel-constraint} is ``imposed" explicitly. This is done by assigning the latest known solution for ${\bm u}_f$ to ${\bm u}_s$ in the solid domain. This projection is typically done using discrete delta functions \cite{Peskin02}. In step 3, the solid domain momentum correction Equation \ref{eq:solid-mom} is used to explicitly calculate ${\bm \lambda}_b$ using the latest known solution. The solid domain boundary condition in Equation \ref{eq:solid-mom-BC} is used to explicitly calculate $\bm \Lambda_s$. In practice, since $\bm \Lambda_s$ is smeared as a body force, its calculation is not a separate step; it is part of the momentum correction Equation \ref{eq:solid-mom}. Often in step 3 approximations are used in the calculation of ${\bm \lambda}$ --- the most common being to neglect ${\bm \sigma}_f$. \REVIEW{The calculation of $\bm \lambda$ can also be simplified by assuming a viscoelastic solid rather than a pure elastic one. Based on this assumption, ${\bm \sigma}_s = {\bm \sigma}_e + {\bm \sigma}_f $, where ${\bm \sigma}_e$ represents the pure elastic component of stress. In both cases, ${\bm \sigma}_f$ is excluded when calculating $\bm \lambda$.} The latest calculated field for ${\bm \lambda}$ is used in step 1 and the calculation process repeats to the next time-step. 

It is also worth mentioning the immersed finite element method (IFEM) formulation here and its similarity to Peskin's IBM approach. In the classical IFEM approach introduced in Zhang et al.~\cite{zhang2004immersed} and Zhang and Gay~\cite{zhang2007immersed}, the variational spaces $\V{v}_f$ and $\V{v}_s$, and consequently the solution spaces  $\V{u}_f$ and $\V{u}_s$,  are taken to be the same. Doing so, modifies the weak form of the extended domain Equation~\ref{eq:combined-weak} to
\begin{eqnarray}
\label{eq:combined-weak-ifem}
&&{}\int\limits_{\Omega}\rho_f{D {\bm u}_f\over D t}\cdot{\bm v}_f dV
 +\int\limits_{\Omega}{\bm\sigma}_f:{\bm D}({\bm v}_f)dV 
 -\int\limits_{\Omega}q_f{\bm\nabla}\cdot{\bm u}_fdV = -\int\limits_{\Omega_s} \biggl( \rho_s  - \rho_f \biggl) {D {\bm u}_f\over D t} \cdot{\bm v}_f dV  \nonumber \\ 
&&{} -\int\limits_{\Omega_s}\Delta {\bm\sigma}:{\bm D}({\bm v}_f)dV  +\int\limits_{\Omega_s}({\bm f}_b + \Delta \rho {\bm g} ) \cdot {\bm v}_f dV
+\int\limits_{\partial\Omega_s}{\bm v}_f \cdot {\bm F}_s dS  \quad \forall\; {\bm v}_f\in\ S_{ev},\;\forall \; q_f\in L^2(\Omega).
\end{eqnarray}
The right hand side of  Equation~\ref{eq:combined-weak-ifem} is typically treated explicitly in the IFEM approach, although semi-implicit~\cite{wang2012semi,wang2019theoretical} and fully-implicit~\cite{newren2007unconditionally,wang2009computational} treatments of the solid domain terms are also possible. Notice that the right hand side terms of  Equation~\ref{eq:combined-weak-ifem} are collectively denoted by ${\bm \lambda}$ in step 3 of the explicit IBM Algorithm~\ref{IBM-algo}. Therefore, Algorithm~\ref{IBM-algo} also describes the explicit IFEM algorithm. The benefit of using separate variational spaces $\V{v}_f$ and $\V{v}_s$ in deriving the immersed formulation is that it helps bring forth the strong form of the extended domain equation (Equation~\ref{eq:solid-mom}). Moreover, the formulation described in Section~\ref{sec_unified_formulation} also allows for the possibility that the fictitious fluid in the solid domain can have a velocity different from the actual solid velocity, as opposed to assuming the same velocity field in both domains as in IFEM.       


\subsection{Fictitious domain method: rigid and self-propelling bodies}  \label{sec_fdm}

\begin{figure}[]
\centering
\subfigure[A two-body system of self-propelling rigid spheres]{
\includegraphics[scale = 0.35]{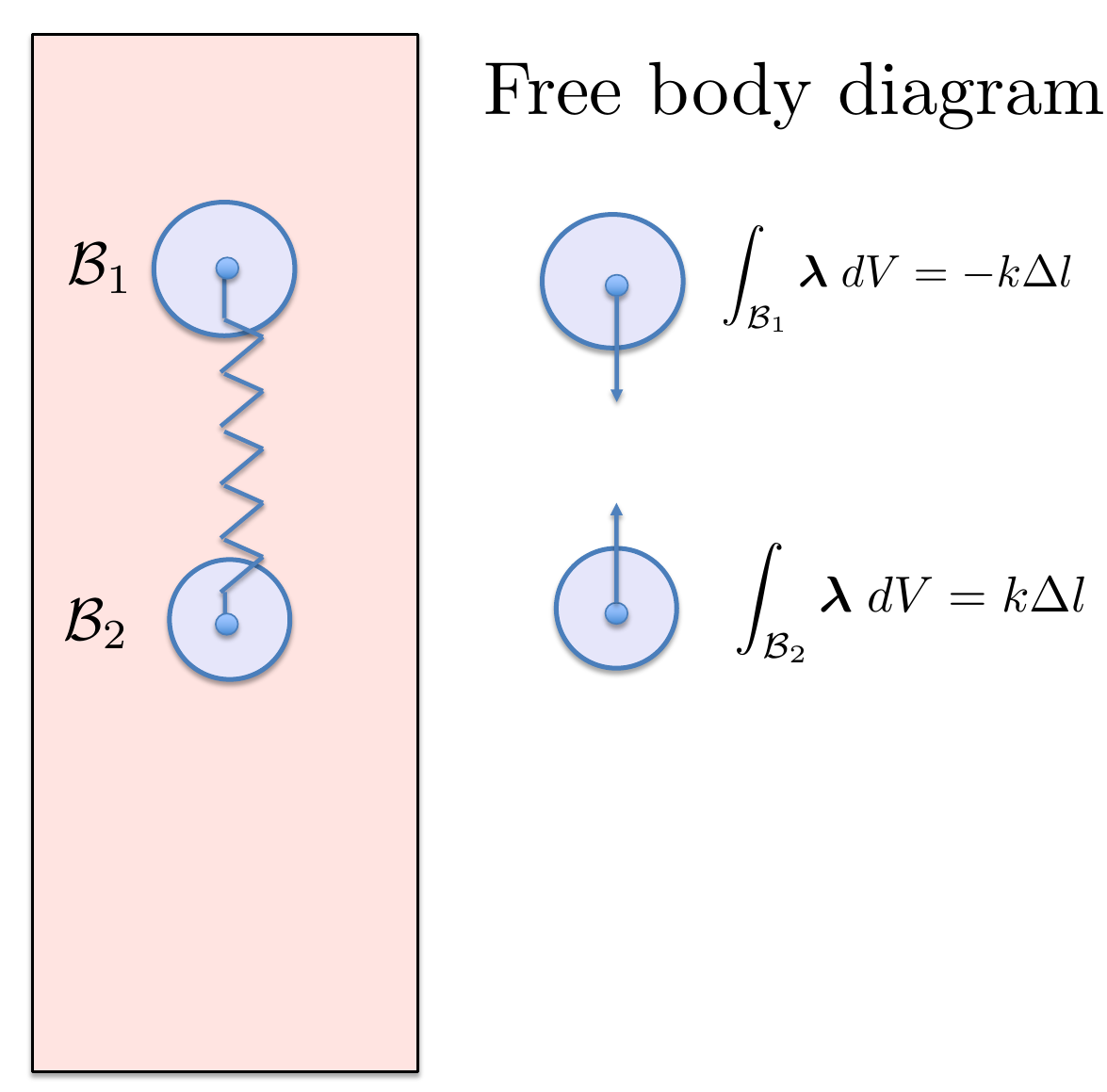}
\label{fig_Two_spheres}
}
 \subfigure[A self-propelling swimmer]{
\includegraphics[scale = 0.35]{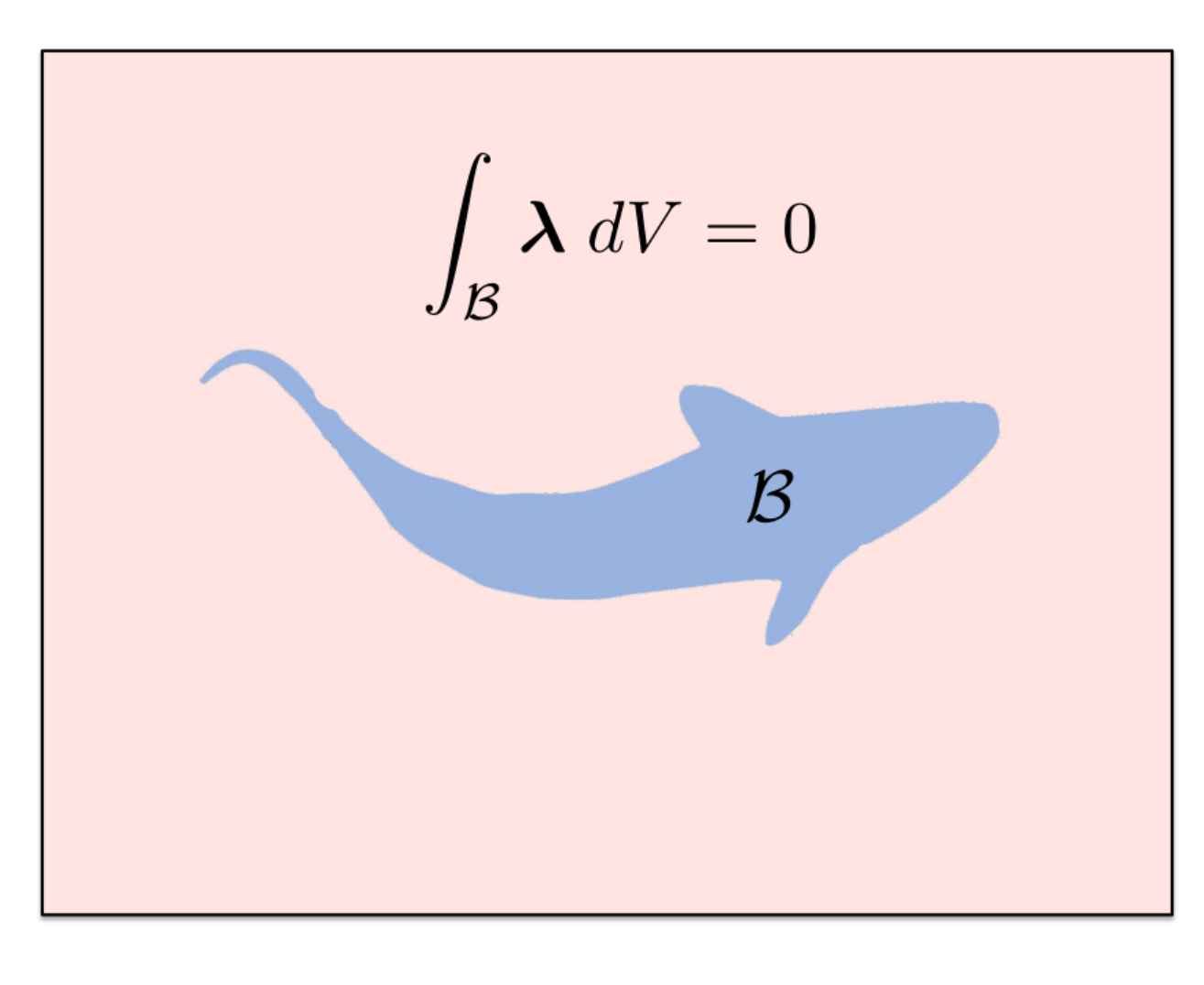}
\label{fig_Fish}
}
  \caption{Representative scenarios of freely moving bodies in fluids:~\subref{fig_Two_spheres} a two-body system of rigid spheres that self-propels itself by actuating a spring tethered at their center of mass points; and~\subref{fig_Fish} a freely-swimming fish that locomotes in water by undulating part of its body.}
\label{fig_FDM}
\end{figure} 

Fictitious domain methods and its variants have been widely used, through the past couple of decades, to simulate freely moving rigid body motion in fluids \cite{glowinski1994fictitious,glowinski1999distributed,Patankar2000,Sharma2005} and self-propelling swimming bodies \cite{Shirgaonkar09,curet2010versatile,APSBhalla13-constraint_ib,bhalla2013forced,patel2018new,bale2014separability,bale2014gray,neveln2014undulating,bhalla2014fully}. Fig.~\ref{fig_FDM} shows two representative scenarios of freely moving bodies in fluids:~\subref{fig_Two_spheres} a two-body system of rigid spheres that self-propels itself by actuating a spring tethered at their center of mass points; and~\subref{fig_Fish} a freely-swimming fish that locomotes in water by undulating part of its body. Here we will consider a self-propelling swimming body whose deformation kinematics are known $\mathit{a\;priori}$, but the translational and angular velocities are not known. Freely moving rigid body motion is a special case where deformation kinematics are zero -- hence that case will not be considered separately. 

Before proceeding further, we need the equations of motion of the solid body (rigid or deforming) as a whole. The translational momentum equation is obtained by integrating the solid momentum correction Equation \ref{eq:solid-mom} and using the solid domain boundary condition in Equation \ref{eq:solid-mom-BC} from the body force formulation. The angular momentum equation is obtained by integrating the moment of the linear momentum Equation \ref{eq:solid-mom}. The resulting equations are
\begin{eqnarray}
\label{eq:solid-linear-mom}
&&{}M_s{D {\bm U}_s\over D t}-{D \over D t}\biggl(\int \limits_{\Omega_s} \rho_f {\bm u}_f dV\biggl) = 
\Delta M {\bm g} - \int \limits_{\mathcal{B}} {\bm \lambda} dV + \int \limits_{\partial\Omega_s} {\bm F}_s dS
+ \int \limits_{\Omega_s} {\bm f}_b dV, \\
\label{eq:solid-angular-mom}
&&{}{D ({\bm I}_s \cdot {\bm \omega}_s) \over D t}-{D \over D t}\biggl(\int \limits_{\Omega_s} ({\bm r} \times \rho_f {\bm u}_f) dV\biggl) = 
- \int \limits_{\mathcal{B}} ({\bm r} \times {\bm \lambda}) dV + \int \limits_{\partial\Omega_s} ({\bm r} \times {\bm F}_s) dS
+ \int \limits_{\Omega_s} ({\bm r} \times {\bm f}_b) dV,
\end{eqnarray}
where $M_s$ and ${\bm I}_s$ are the mass and moment of inertia, respectively, of the solid body, $\Delta M = \Delta \rho V_b$, $V_b$ is the volume of the solid body, and ${\bm r}$ is the coordinate with respect to the center of mass of the solid body. The average translational velocity, ${\bm U}_s$, and the average angular velocity, ${\bm \omega}_s$, of the solid body are defined by
\begin{eqnarray}
\label{eq:trans-vel-defn}
&&{}M_s{\bm U}_s=\int \limits_{\Omega_s} \rho_s {\bm u}_s dV, \\
\label{eq:ang-vel-defn}
&&{}{\bm I}_s \cdot {\bm \omega}_s = \int \limits_{\Omega_s} ({\bm r} \times \rho_s {\bm u}_s) dV.
\end{eqnarray}
Note the notation pertaining to integrals involving $\bm \lambda$
\begin{eqnarray}
\label{eq:lambda-linear-int}
&&{}\int \limits_{\mathcal{B}} {\bm \lambda} dV = \int \limits_{\partial\Omega_s} {\bm \Lambda}_s dS 
+ \int \limits_{\Omega_s} {\bm \lambda}_b dV, \\
\label{eq:lambda-ang-int}
&&{}\int \limits_{\mathcal{B}} ({\bm r} \times {\bm \lambda}) dV =  \int \limits_{\partial\Omega_s} ({\bm r} \times {\bm \Lambda}_s) dS 
+  \int \limits_{\Omega_s} ({\bm r} \times {\bm \lambda}_b) dV,
\end{eqnarray}
When the constraint ${\bm u}_f = {\bm u}_s$ on $\mathcal{B}(t)$ (Equation \ref{eq:vel-constraint}) is used, it is convenient to write Equations \ref{eq:solid-linear-mom} and \ref{eq:solid-angular-mom}, after replacing ${\bm u}_f$ with ${\bm u}_s$, as follows
\begin{eqnarray}
\label{eq:solid-linear-mom-inertia}
&&{} \Delta M_s{D {\bm U}_s\over D t} = - \int \limits_{\mathcal{B}} {\bm \lambda} dV + \int \limits_{\partial\Omega_s} {\bm F}_s dS
+ \int \limits_{\Omega_s} {\bm f}_b dV + \Delta M {\bm g}, \\
\label{eq:solid-ang-mom-inertia}
&&{}{D (\Delta {\bm I}_s \cdot {\bm \omega}_s) \over D t} = - \int \limits_{\mathcal{B}} ({\bm r} \times {\bm \lambda}) dV + 
\int \limits_{\partial\Omega_s} ({\bm r} \times {\bm F}_s) dS + \int \limits_{\Omega_s} ({\bm r} \times {\bm f}_b) dV,
\end{eqnarray}
where $\Delta {\bm I}$ is the moment of inertia of the body based on the density difference $\Delta \rho$. Note that for a neutrally buoyant body (which is a reasonable assumption for swimming animals) the inertia terms, proportional to $\Delta \rho$, are zero. The fictitious domain method (FDM) for a neutrally buoyant self-propelling body is summarized in Algorithm~\ref{FDM-algo}.
  \begin{algorithm}[]
  \caption{FICTITIOUS DOMAIN METHOD: SELF-PROPULSION} \label{FDM-algo}
  \begin{algorithmic}
    \State \textbf{1.} Solve  \\ \\
    $\quad \quad \rho_f {D {\bm u}_f\over D t}= {\bm\nabla}\cdot{\bm\sigma}_f+{\bm \lambda}~\textrm{in}~\Omega$,  \\ \\
    $\quad \quad {\bm u}_f  = {\bm U}_s + ({\bm \omega}_s \times {\bm r}) + {\bm u}_{s,\textrm{def}}~\textrm{on}~\mathcal{B}(t)$,  
    \Comment{ ${\bm u}_{s,\textrm{def}}$ are the known deformation kinematics} \\  \\
    
    $ \quad \quad \int \limits_{\mathcal{B}} {\bm \lambda} dV = \int \limits_{\partial\Omega_s} {\bm F}_s dS + \int \limits_{\Omega_s} {\bm f}_b dV $,
      \Comment{$ \quad \quad \int \limits_{\mathcal{B}} {\bm \lambda} dV = {\bm 0}$ in the absence of external forces}   \\ \\

  $\quad \quad \int \limits_{\mathcal{B}} ({\bm r} \times {\bm \lambda}) dV = 
\int \limits_{\partial\Omega_s} ({\bm r} \times {\bm F}_s) dS + \int \limits_{\Omega_s} ({\bm r} \times {\bm f}_b) dV$.  
\Comment{$\int \limits_{\mathcal{B}} ({\bm r} \times {\bm \lambda}) dV = {\bm 0}$ in the absence of external torques}  \\

   \State \textbf{2.} Compute \\ \\
   $\quad \quad {\bm f_b} = \Delta\rho{D {\bm u}_s\over D t} - {\bm\nabla}\cdot({{\bm\sigma}_s - {\bm\sigma}_f})+{\bm \lambda_b} - \Delta\rho{\bm g}~\textrm{in}~\Omega_s$, \\ \\
   $\quad \quad {\bm F}_s = \biggl[{\bm\sigma}_s^- - {\bm \sigma}_f^-\biggl]\cdot\hat{\bm n} + {\bm \Lambda}_s 
~\textrm{on}~\partial\Omega_s$. \Comment{$\f_b$ and $\F_s$ are postprocessed quantities}

  \end{algorithmic}
  \end{algorithm}

In FDM for self-propulsion, the fluid momentum Equation \ref{eq:fluid-mom} of the body force formulation is solved coupled with the velocity constraint Equation \ref{eq:vel-constraint}. The solid velocity field is of the form ${\bm u}_s = {\bm U}_s + ({\bm \omega}_s \times {\bm r}) + {\bm u}_{s,\textrm{def}}$, where  ${\bm u}_{s,\textrm{def}}$ is known. Hence, only the rigid component of the overall solid velocity field, i.e., ${\bm U}_s$ and ${\bm \omega}_s$ are unknown. Two more equations are needed to solve for ${\bm U}_s$ and ${\bm \omega}_s$. These are obtained from Equations \ref{eq:solid-linear-mom-inertia} and \ref{eq:solid-ang-mom-inertia}, which simplify for a neutrally buoyant body to the form seen in equation set of Algorithm~\ref{FDM-algo} (equations involving integrals of $\bm \lambda$). Note that Equations \ref{eq:solid-linear-mom-inertia} and \ref{eq:solid-ang-mom-inertia} were obtained after integrating the solid momentum correction Equation \ref{eq:solid-mom} together with the solid boundary condition Equation \ref{eq:solid-mom-BC} in the body force formulation. The equations involving integral of $\bm \lambda$ in equation set of Algorithm~\ref{FDM-algo} are further simplified by noting that the integral of the external forces and moments on a freely swimming body are zero as shown in Fig.~\ref{fig_Fish}.  Hence, the integrals of ${\bm F}_s$ and ${\bm f}_b$, and the integrals of their moments are zero. This gives relations that the integral of $\bm \lambda$ and the integral of the moment of $\bm \lambda$ are both zero. In contrast, for the self-propelling two sphere system shown in Fig.~\ref{fig_Two_spheres}, the integral of $\bm \lambda$ for each sphere is non-zero and is equal to the net actuator (spring) force. These additional requirements provide two equations needed to solve for ${\bm U}_s$ and ${\bm \omega}_s$ at each instant. Finally, once the solution is obtained, Equations \ref{eq:solid-mom} and \ref{eq:solid-mom-BC} can be applied pointwise to calculate $\mathit{a\;posteriori}$ the force fields ${\bm f}_b$ and ${\bm F}_s$ if the elastic properties of the body are known. This is the postprocess step 2 in the equation set of Algorithm~\ref{FDM-algo} above. These force fields provide insight into what force field might have been generated by the body to propel itself according to the specified deformation kinematics.

Having seen how the FDM follows from the body force formulation in Equations \ref{eq:fluid-mom}--\ref{eq:solid-mom-BC} we briefly consider how these equations are solved. In practice the FDM algorithm is typically executed as demonstrated in Algorithm~\ref{FDM-Chorin-algo} using a first-order fractional time-stepping  (FTS) scheme.  The algorithm proceeds by solving the fluid equations in the entire domain as shown in step 1. Step 2 essentially solves the solid momentum correction equation (Equations \ref{eq:solid-mom} and \ref{eq:solid-mom-BC} in the body force formulation) in an integrated form to get  ${\bm U}_s$, ${\bm \omega}_s$, and $\bm \lambda$. In step 3, the imposition of the velocity constraint is executed by simply replacing $\widehat{\bm u}_f$ field with latest known ${\bm u}_s$ field from step 2. This is equivalent to adding a force field $\bm \lambda$ in the fluid. Similar to implementation of IBM, $\bm \Lambda_s$ is typically smeared over few grids at the fluid-solid interface using regularized integral kernels \cite{Peskin02}. Thus, here too, in practice the entire $\bm \lambda$ field is treated like a body force term.

  \begin{algorithm}[]
  \caption{FICTITIOUS DOMAIN METHOD: SELF-PROPULSION FTS 
  } \label{FDM-Chorin-algo}
  \begin{algorithmic}
    \State \textbf{1.} Solve fluid momentum equation in the entire domain $\Omega$  \\ \\
    $\quad \quad \rho_f { {\widehat{\bm u}_f - {\bm u}_f^n}\over \Delta t} + \rho_f( \widehat{\bm u}_f \cdot \bm\nabla)\widehat{\bm u}_f = {\bm\nabla}\cdot{\widehat{\bm \sigma}}_f$, \\ \\
    $\quad \quad \bm\nabla \cdot \widehat{\bm u}_f = \V{0}$.
    \Comment{\; $\widehat{}$ denotes intermediate fields and superscript $n$ is the time-step counter} \\
     
    \State \textbf{2.} Solve the rigid velocity component of the solid \\ \\
    $\quad \quad \int \limits_{\mathcal{B}} {\bm \lambda}^{n+1} dV = {\bm 0} \implies {\bm U}_s^{n+1}$, 
    \Comment{external forces are assumed to be absent} \\
    
    $\quad \quad	\int \limits_{\mathcal{B}} ({\bm r} \times {\bm \lambda}^{n+1}) dV = {\bm 0} \implies {\bm \omega}_s^{n+1}$, 
    \Comment{external torques are assumed to be absent} \\
    
    $\quad \quad {\bm \lambda}^{n+1} = \rho_f { {{\bm U}_s^{n+1} + ({\bm \omega}_s^{n+1} \times {\bm r}) + {\bm u}_{s,\textrm{def}}^{n+1} - \widehat{\bm u}_f }\over \Delta t}$. \Comment{$\V{\lambda}$ defined on $\mathcal{B}(t)$}.  \\
    
    \State \textbf{3.} Impose velocity constraint on $\mathcal{B}(t)$ \\  \\
    $\quad \quad \rho_f { {{\bm u}_f^{n+1} - \widehat{\bm u}_f }\over \Delta t} = {\bm \lambda}^{n+1}.$  
     \Comment{Chorin-type projection}
    
  \end{algorithmic}
  \end{algorithm}


\subsection{Velocity/direct forcing method}  \label{sec_vfm}

Velocity forcing methods are those in which the solid domain velocity field is directly assigned to the fluid either on the fluid-solid interface or in the entire solid domain. This technique has been used to perform simulations where the solid velocity is known $\mathit{a\;priori}$ \cite{fadlun2000combined} or to solve for the motion of freely-moving rigid bodies in fluids \cite{uhlmann2005immersed}, among others. To explain this method in the context of the body force formulation (Equations \ref{eq:fluid-mom}--\ref{eq:solid-mom-BC}) we will consider a freely-moving rigid body in a fluid. The algorithm (Algorithm~\ref{VFM-algo}) for the velocity forcing method proceeds as follows. In step 1,  the fluid momentum Equation \ref{eq:fluid-mom} is solved coupled with the velocity constraint Equation \ref{eq:vel-constraint}. The latest known solid velocity field ${\bm u}_s$ is used in the solid domain to impose the constraint. The constraint can be imposed on $\partial\Omega_s(t)$ or on $\mathcal{B}(t)$. Note that the solid velocity is of the form ${\bm u}_s = {\bm U}_s + {\bm \omega}_s \times {\bm r}$ since we are considering a rigid body as an example. In step 2 of the velocity forcing method, the translational and angular velocities of the solid body (Equations \ref{eq:solid-linear-mom} and \ref{eq:solid-angular-mom}) are explicitly calculated using the latest known fluid velocity field ${\bm u}_f$ and the Lagrange multiplier field $\bm \lambda$. The execution of step 2 is equivalent to solving the integrated form of Equations \ref{eq:solid-mom} and \ref{eq:solid-mom-BC} in the body force formulation. We remark that the net hydrodynamic force and torque acting on the rigid body can also be evaluated via integral of pointwise hydrodynamic force and torque acting on the surface of the solid from the fluid side. This approach is also frequently used in the literature. 

 \begin{algorithm}[]
  \caption{VELOCITY FORCING METHOD} \label{VFM-algo}
  \begin{algorithmic}
    \State \textbf{1.} Solve fluid momentum equation coupled with the velocity constraint  \\ \\
	$\quad \quad \rho_f{D {\bm u}_f\over D t} = {\bm\nabla}\cdot{\bm\sigma}_f+{\bm \lambda}~\textrm{in}~\Omega$, \\ \\
	$\quad \quad {\bm u}_f = {\bm U}_s + {\bm \omega}_s \times {\bm r}~\textrm{on}~\partial\Omega_s(t)~\textrm{or on}~\mathcal{B}(t)$. \\ 
    
    \State \textbf{2.} Update the rigid body velocities  \\ \\
	$\quad \quad M_s{D {\bm U}_s\over D t} = \boxed{{D \over Dt} \left(\;\int \limits_{\Omega_s} \rho_f {\bm u}_f dV\right)  - \int \limits_{\Omega_s} {\bm \lambda} dV}  + \int \limits_{\partial\Omega_s} {\bm F}_s dS + \int \limits_{\Omega_s} {\bm f}_b dV + \Delta M {\bm g}$,  \\
	
$\quad \quad {D ( {\bm I}_s \cdot {\bm \omega}_s) \over D t} = \boxed{{D \over Dt} \left(\; \int \limits_{\Omega_s}{\bm r} \times \rho_f {\bm u}_f dV\right) - \int \limits_{\Omega_s} ({\bm r} \times {\bm \lambda}) dV}  + \int \limits_{\partial\Omega_s} ({\bm r} \times {\bm F}_s) dS + \int \limits_{\Omega_s} ({\bm r} \times {\bm f}_b) dV$,  \\ \\
which follow from Equations~\ref{eq:solid-linear-mom} and~\ref{eq:solid-angular-mom}. Here the boxed terms are the net hydrodynamic force and torque acting on the body. This can be shown by using the Gauss-divergence theorem on Equation~\ref{eq:fluid-mom} as \\ \\
$\quad \quad {D \over Dt} \left(\; \int \limits_{\Omega_s} \rho_f {\bm u}_f dV\right)  - \int \limits_{\Omega_s} {\bm \lambda} dV =   \int \limits_{\partial\Omega_s} {\bm \sigma}_f^{+} \cdot \V{\hat{n}} dS $, \\

$ \quad \quad {D \over Dt} \left(\; \int \limits_{\Omega_s}{\bm r} \times \rho_f {\bm u}_f dV\right) - \int \limits_{\Omega_s} ({\bm r} \times {\bm \lambda}) dV = \int \limits_{\partial\Omega_s} \V{r} \times ({\bm \sigma}_f^{+} \cdot \V{\hat{n}}) dS$. \\ \\

\State \textbf{3.} Go to step 1.   
    
  \end{algorithmic}
  \end{algorithm}

\begin{algorithm}[]
  \caption{VELOCITY FORCING METHOD: STEP 1 IMPLEMENTATION} \label{VFM-step1-algo}
  \begin{algorithmic}
    \State \textbf{1(a)} Define the Lagrange multiplier field \\ \\
        $\quad \quad {\bm \lambda}^{n+1} =  \kappa  \left( {\bm U}_s + ({\bm \omega}_s \times {\bm r}) - \widehat{\bm u}_f  \right)~\textrm{on}~\partial\Omega_s~\textrm{or}~\mathcal{B}(t)$. \Comment{$\kappa \approx \rho_f/\Delta t$}  \\

    \State \textbf{1(b)} Solve fluid momentum equation in the entire domain $\Omega$  \\ \\
    
$\quad \quad \rho_f { {\widehat{\bm u}_f - {\bm u}_f^n}\over \Delta t} + \rho_f( \widehat{\bm u}_f \cdot \bm\nabla)\widehat{\bm u}_f = {\bm\nabla}\cdot{\widehat{\bm \sigma}}_f + \theta {\bm \lambda}^{n+1}$,  
\Comment{$\theta = 1$ for BP method and $\theta = 0$ for FTS scheme}
\\  \\
 $\quad \quad \bm\nabla \cdot \widehat{\bm u}_f = 0$.  \\
    
    \State \textbf{1(c)} \textbf{if} ($\theta = 0$) \textbf{then} \\ \\
     $\quad \quad \quad \rho_f { {{\bm u}_f^{n+1} - \widehat{\bm u}_f }\over \Delta t} = {\bm \lambda}^{n+1}~\textrm{on}~\partial\Omega_s(t)~\textrm{or}~\mathcal{B}(t)$. \Comment{Chorin-type projection}  \\ \\
       $\quad \quad$\textbf{endif} 
       
  \end{algorithmic}
  \end{algorithm}

Step 1 of the velocity forcing algorithm can be implemented either as a first-order fractional time-stepping (FTS) scheme or as a Brinkman penalization (BP) method~\cite{angot1999penalization,bhalla2020simulating,thirumalaisamy2023effective} --- as described by Algorithm~\ref{VFM-step1-algo}.  Briefly,  the Lagrange multiplier field $\V{\lambda}$ is defined using the latest rigid body velocities ${\bm U}_s$ and ${\bm \omega}_s$ as noted in step 1(a). In step 1(b), the momentum equation is solved in the entire domain. Step 1(c) is required only for the FTS scheme, which imposes the velocity constraint by simply replacing $\widehat{\bm u}_f$ field with the latest known ${\bm u}_s$ field at the location of the constraint. This is equivalent to adding a force field $\bm \lambda$ in the fluid, as done for the BP implementation in step 1(b). In the literature, methods which directly replace $\widehat{\bm u}_f$ with ${\bm u}_s$ using a FTS scheme are often times referred to as \emph{direct forcing} methods~\cite{uhlmann2005immersed,yang2012simple} instead of velocity forcing methods.

\subsection{Fully implicit method and Brownian simulations}  \label{sec_stokes}

In Section~\ref{sec_fdm} we described the fictitious domain method for a self-propelling body, and a typical fractional time-stepping (FTS) scheme to implement it. Although the FTS approach allows for an efficient implementation of the fictitious domain method, the constraint in the extended fluid region is imposed only approximately. Moreover, the FTS approach does not work in the limit of zero Reynolds number (steady Stokes regime), as the inertial terms in Equations~\ref{eq:fluid-mom} and~\ref{eq:solid-mom} are absent. These two limitations motivate solving the fictitious domain algorithm implicitly. Doing so requires solving for the fluid velocity $\u_f$, pressure $p_f$, Lagrange multiplier field $\V{\lambda}$, and the rigid body velocities $\V{U}_s$ and $\V{\omega}_s$ (we consider a rigid body in this discussion), together as a simultaneous system of equations. Algorithm~\ref{PresImplicit-FDM-algo} describes the implicit fictitious domain method, when the rigid body translational and rotational velocity is prescribed (i.e. is known), and Algorithm~\ref{SelfImplicit-FDM-algo} describes the implicit fictitious domain method when the rigid body velocities are also sought in the solution.

  \begin{algorithm}[]
  \caption{FULLY IMPLICIT FICTITIOUS DOMAIN METHOD: PRESCRIBED VELOCITY} \label{PresImplicit-FDM-algo}
  \begin{algorithmic}
    \State \textbf{1.} Solve  the three simultaneous system of equations for $\u_f$, $p_f$, and $\V{\lambda}$ \\ \\
    $\quad \quad \rho_f {D {\bm u}_f\over D t}= {\bm\nabla}\cdot{\bm\sigma}_f+{\bm \lambda}~\textrm{in}~\Omega$,  
    \Comment{$ \rho_f {D {\bm u}_f\over D t}$ term is omitted in the steady Stokes limit}\\ \\
    $\quad \quad {\bm\nabla}\cdot \u_f = \V 0$, \Comment {$p_f$ is the Lagrange multiplier  enforcing incompressibility constraint} \\ \\
    $\quad \quad {\bm u}_f   = {\bm U}_s + ({\bm \omega}_s \times {\bm r}) + {\bm u}_{s,\textrm{def}}~\textrm{on}~\mathcal{B}(t)$.  \Comment{$\V{\lambda}$ is the Lagrange multiplier enforcing $\u_f = \u_s$ on $\mathcal{B}(t)$} \\
    \Comment{ $\U_s, \V{\omega}_s$ and ${\bm u}_{s,\textrm{def}}$ are known $\mathit{a\;priori}$} 
   
  \end{algorithmic}
  \end{algorithm}
  
    \begin{algorithm}[]
  \caption{FULLY IMPLICIT FICTITIOUS DOMAIN METHOD: SELF-PROPULSION} \label{SelfImplicit-FDM-algo}
  \begin{algorithmic}
    \State \textbf{1.} Solve  the five simultaneous system of equations for $\u_f$, $p_f$, $\V{\lambda}$, $\V{U}_s$ and $\V{\omega}_s$ \\ \\
    $\quad \quad \rho_f {D {\bm u}_f\over D t}= {\bm\nabla}\cdot{\bm\sigma}_f+{\bm \lambda}~\textrm{in}~\Omega$,  
     \Comment{$ \rho_f {D {\bm u}_f\over D t}$ term is omitted in the steady Stokes limit}\\ \\
    $\quad \quad {\bm\nabla}\cdot \u_f = \V 0$, \Comment {$p_f$ is the Lagrange multiplier enforcing incompressibility constraint} \\ \\
    $\quad \quad {\bm u}_f   = {\bm U}_s + ({\bm \omega}_s \times {\bm r}) + {\bm u}_{s,\textrm{def}}~\textrm{on}~\mathcal{B}(t)$,  \Comment{$\V{\lambda}$ is the Lagrange multiplier enforcing $\u_f = \u_s$ on $\mathcal{B}(t)$} \\ \\
        
    $ \quad \quad \int \limits_{\mathcal{B}} {\bm \lambda} dV = \int \limits_{\partial\Omega_s} {\bm F}_s dS + \int \limits_{\Omega_s} {\bm f}_b dV $,
      \Comment{constraint equation for translational velocity $\U_s$}   \\ \\

  $\quad \quad \int \limits_{\mathcal{B}} ({\bm r} \times {\bm \lambda}) dV = 
\int \limits_{\partial\Omega_s} ({\bm r} \times {\bm F}_s) dS + \int \limits_{\Omega_s} ({\bm r} \times {\bm f}_b) dV$.  
\Comment{constraint equation for rotational velocity $\V{\omega}_s$}

  \end{algorithmic}
  \end{algorithm}
  
  The fully implicit system described in Algorithms~\ref{PresImplicit-FDM-algo} and~\ref{SelfImplicit-FDM-algo} can be solved using specialized physics-based preconditioners~\cite{kallemov2016immersed,balboa2017hydrodynamics} or using iterative solvers like SIMPLE/PISO~\cite{curet2010versatile}. \REVIEW{Typically, fully implicit IB methods are used to simulate the Brownian motion of particles in overdamped (Stokes) flow regimes~\cite{sharma2004direct, Chen2006Brownian, sprinkle2019brownian}; see next paragraph on some technical issues related to modeling Brownian motion of particles. The majority of IB applications reported in the literature use (efficient) FTS approaches to deal with flows at moderate to high Reynolds numbers. There are a few IB methods reported for zero Reynolds numbers in the literature~\cite{kallemov2016immersed,balboa2017hydrodynamics}. Solving the resulting large system of equations is also challenging in this case. It should be noted, however, that the body force formulation (Equations \ref{eq:fluid-mom}--\ref{eq:solid-mom-BC}) is Reynolds number independent.
  
 The distributed Lagrange multiplier method allow natural extension to simulate Brownian motion of particles. Assuming the entire fluid domain to be a fluctuating fluid and constraining the particle domains to move rigidly, leads to Brownian motion of the immersed particles \cite{patankar2001brownian, sharma2004direct}. Short time scale Brownian calculations naturally capture the correct algebraic tails in the translation and rotational autocorrelation functions as opposed to the exponential tail in the solution of the Langevin equation for Brownian motion \cite{Chen2006Brownian}. Yet, few aspects need careful consideration. First, to satisfy the fluctuation-dissipation theorem for the discretized fluid equations, it is essential to use a spatial discretization scheme with good spectral properties. Hence, a staggered discretization scheme for the fluid variables is recommended \cite{sharma2004direct, Chen2006Brownian}. Second, to get the correct statistics of particle diffusion on the long time scale, special integrators are essential to correctly capture the effect of the changes in particle configuration \cite{sprinkle2019brownian}.}



\section{Multiphase formulations} \label{sec_multiphase_system}

\begin{figure}[]
\centering
\subfigure[Solid-aware gas-liquid advection]{
\includegraphics[scale = 0.42]{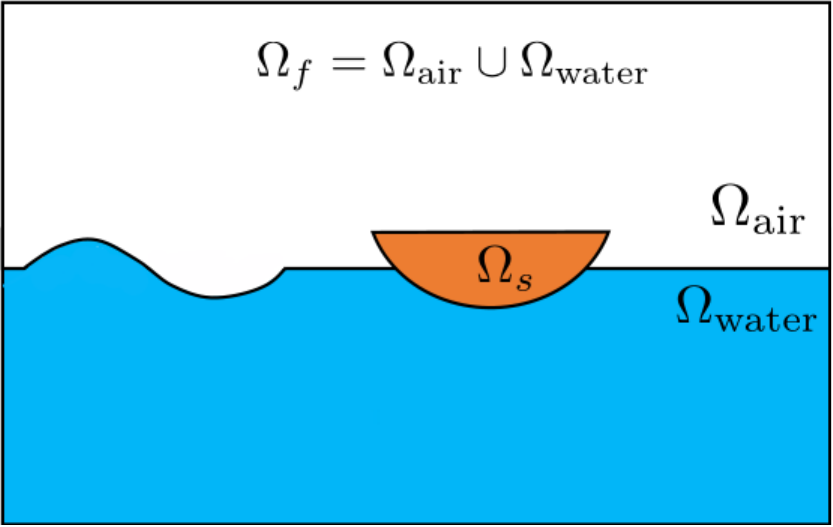}
\label{fig_wec_vof}
}
 \subfigure[Solid-agnostic gas-liquid advection]{
\includegraphics[scale = 0.42]{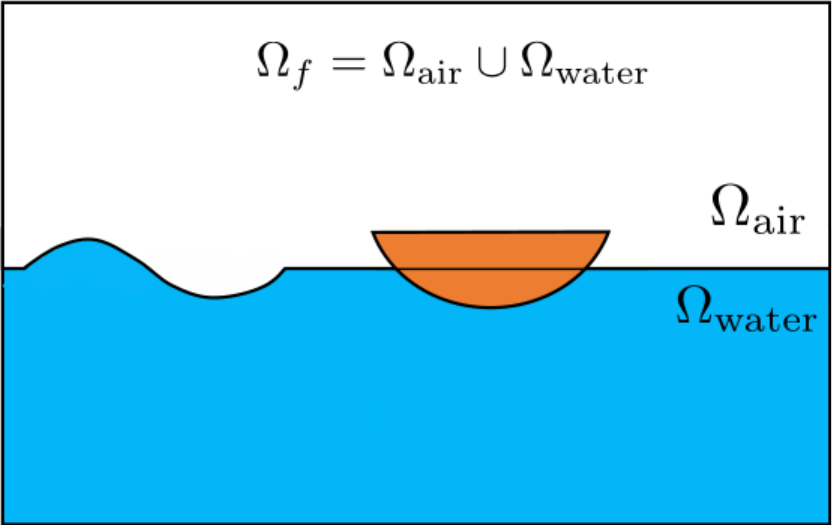}
\label{fig_wec_ls}
}
  \caption{A schematic representation of a wave energy converter device oscillating on the air-water interface under the action of incoming water waves. The air phase $\Omega_{\rm air}$ is shown in white, the water phase $\Omega_{\rm water}$ in blue, and the solid phase $\Omega_{\rm solid}$ in brown color, respectively. The flowing/fluid phase domain is $\Omega_f(t) = \Omega_{\rm air}(t) \cup \Omega_{\rm water}(t)$. The density ratio between air and water phase is approximately 1000, whereas between air and solid (mechanical oscillator) phase is approximately 500. The gas-liquid interface is advected in \subref{fig_wec_vof}  $\Omega_f(t)$  and \subref{fig_wec_ls} $\Omega$ regions, respectively. }
\label{fig_wec}
\end{figure} 


The extended domain system discussed in Section~\ref{sec_unified_formulation} assumed a constant value of fluid and solid density $\rho_f$ and $\rho_s$, respectively. We did so for clarity of exposition and because the majority of prior literature made that assumption. This assumption can be relaxed in both the weak and strong form of the equations derived in Section~\ref{sec_unified_formulation} to enable simultaneous modeling of solid, liquid and gas phases. With $\rho_f(\x,t)$ denoting the density of ``flowing" phases (such as air and water) in $\Omega_f(t)$, and $\rho_s(\x,t)$ denoting the density of the solid in $\Omega_s(t)$, the strong form of the multiphase system reads as
 \begin{eqnarray}
\label{eq:multi-eqns1}
&&{}\rho_f(\x,t){D {\bm u}_f\over D t} = {\bm\nabla}\cdot{\bm\sigma}_f + \rho_f(\x,t)\V{g}~\textrm{in}~\Omega_f(t),\\
\label{eq:multi-eqns2}
&&{}\rho_s(\x,t){D {\bm u}_s\over D t} = {\bm\nabla}\cdot{\bm\sigma}_s +  {\bm f}_b +  \rho_s(\x,t) {\bm g}~\textrm{in}~\Omega_s(t),\\
\label{eq:multi-eqns3}
&&{}{\bm\nabla}\cdot{\bm u}_f=0~\textrm{in}~\Omega_f(t), {\bm\nabla}\cdot{\bm u}_s=0~\textrm{in}~\Omega_s(t),\\
\label{eq:multi-eqns4}
&&{}{D\rho_f (\x,t)\over D t} = 0~\textrm{in}~\Omega_f(t),{D\rho_s(\x,t)\over D t} = 0~\textrm{in}~\Omega_s(t), \\
\label{eq:multi-eqns5}
&&{}{\bm u}_f = {\bm u}_s~\textrm{on}~\partial\Omega_s(t),\\
\label{eq:multi-eqns6}
&&{}\biggl[{\bm\sigma}_f^+ - {\bm \sigma}_s^-\biggl]\cdot\hat{\bm n} =- {\bm F}_s~\textrm{on}~\partial\Omega_s(t).
\end{eqnarray}
We have included the gravitational force separately in the flowing/fluid phase region $\Omega_f$ and in the solid phase region $\Omega_s$. This is to account for situations where the body is partially submerged in different fluid phases. This effect can be observed in the case of a wave energy converter device oscillating (heaving and pitching) on the air-water interface, as shown in Fig.~\ref{fig_wec}. In the fluid momentum equation, surface tension could be included as a body force term along with other forces if desired. A volume conservation equation is represented by Equation~\ref{eq:multi-eqns3}. Additionally, Equation~\ref{eq:multi-eqns4}, which represents the incompressibility of the three phases, is needed to track the evolving domain of each phase. Depending on the computational technique, these equations are substituted by different interface advection/tracking equations. \REVIEW{A discussion on the interface tracking methods is provided after describing the computational algorithms for solving Equations~\ref{eq:multi-eqns1}--\ref{eq:multi-eqns6} later in this section}. Based on the steps in Section~\ref{sec_unified_formulation} for a constant density system, the extended domain strong form for the variable density multiphase system is given by
\begin{eqnarray}
\label{eq:multi:fluid-mom}
&&{}\rho_f(\x,t){D {\bm u}_f\over D t} = {\bm\nabla}\cdot{\bm\sigma}_f + \rho_f(\x,t)\V{g}  + {\bm \Lambda_s} \delta_s + {\bm \lambda_b}~\textrm{in}~\Omega, \\
\label{eq:multi:solid-mom}
&&{}\rho_s(\x,t) {D {\bm u}_s\over D t} - \rho_f(\x,t) {D {\bm u}_f \over D t} = {\bm\nabla}\cdot\Delta{\bm\sigma}+{\bm f_b} + (\rho_s(\x,t)  - \rho_f(\x,t)) {\bm g}-{\bm \lambda_b}~\textrm{in}~\Omega_s(t), \\
\label{eq:multi:vel-constraint}
&&{}{\bm u}_f = {\bm u}_s~\textrm{on}~\mathcal{B}(t), \\
\label{eq:multi:solid-mom-BC}
&&{}\Delta\bm\sigma\cdot\hat{\bm n} =
\biggl[{\bm\sigma}_s^- - {\bm \sigma}_f^-\biggl]\cdot\hat{\bm n} = {\bm F}_s - {\bm \Lambda}_s~\textrm{on}~\partial\Omega_s(t).
\end{eqnarray}
Next, we will discuss  the fictitious domain and  velocity forcing algorithms to solve the above equation system.

\subsection{Fictitious domain method: multiphase system}

Let the density of the fluid extended into the solid region (fictitious fluid) be matched to the solid density, i.e., $\rho_f(\x,t) = \rho_s(\x,t)$ in $\mathcal{B}(t)$. With this choice, the fictitious domain method for the variable density multiphase system is the same as Algorithm~\ref{FDM-Chorin-algo} except that the fluid density is, in general, not constant over the entire domain. A variable density fluid solver is therefore required rather than a (possibly faster) constant density solver. The matching of the fictitious fluid and solid densities implies that the inertial correction and gravitational force terms in Equation~\ref{eq:multi:solid-mom} are absent.  \REVIEW{Algorithm~\ref{FDM-Chorin-algo} is used in Pathak and Raessi~\cite{pathak20163d} and Nangia et al.~\cite{nangia2019dlm} to solve multiphase FSI equations.} 

\subsection{Velocity forcing method: multiphase system}

The velocity forcing method for the multiphase FSI system remains the same as described in Algorithms~\ref{VFM-algo} and~\ref{VFM-step1-algo}. However, in this case, it is not required to match the fictitious fluid density with solid density within the solid region $\Omega_s$. As a result, any reasonable value of fictitious fluid density can be used within the solid domain. Note that the first boxed term in the velocity forcing Algorithm~\ref{VFM-algo} contains the contribution of net hydrodynamic and gravitational forces because of the modification introduced in Equation~\ref{eq:multi:fluid-mom}, i.e. in this case 
\begin{equation}
{D \over Dt} \left(\; \int \limits_{\Omega_s} \rho_f {\bm u}_f dV \right)  - \int \limits_{\Omega_s} {\bm \lambda} dV =   \int \limits_{\partial\Omega_s} {\bm \sigma}_f^{+} \cdot \V{\hat{n}} dS +  \int \limits_{\Omega_s} \rho_f(\x,t)\V{g} dV. 
\end{equation}
Recalling the definition of $\Delta M =  \int \limits_{\Omega_s} \left( \rho_s(\x,t) - \rho_f(\x,t) \right) dV$, the rigid body translation equation in Algorithm~\ref{VFM-algo} can be directly written as 
\begin{align}
 M_s{D {\bm U}_s\over D t} &=   \int \limits_{\partial\Omega_s} {\bm \sigma}_f^{+} \cdot \V{\hat{n}} dS   + \int \limits_{\partial\Omega_s} {\bm F}_s dS + \int \limits_{\Omega_s} {\bm f}_b dV + +  \int \limits_{\Omega_s} \rho_f(\x,t)\V{g} dV + \Delta M {\bm g},  \nonumber  \\
 & = \int \limits_{\partial\Omega_s} {\bm \sigma}_f^{+} \cdot \V{\hat{n}} dS + \int \limits_{\partial\Omega_s} {\bm F}_s dS + \int \limits_{\Omega_s} {\bm f}_b dV + M_s \V{g}.
\end{align}
The rigid body angular momentum equation remains the same because the gravitational force term in Equation~\ref{eq:multi:fluid-mom} does not produce any rotational torque on the body. \REVIEW{Algorithms~\ref{VFM-algo} and~\ref{VFM-step1-algo} is used in~\cite{bergmann2022numerical,BhallaBP2019,Bergmann2015,Khedkar2020,khedkar2022model} to simulate multiphase FSI in conjunction with the Brinkman penalization technique.}

\REVIEW{
\subsection{Interface tracking methods}
While there are several approaches for tracking the interface between two fluid phases implicitly, the level set (LS) method and the volume of fluid (VOF) method are two popular choices. Within the context of three phase gas-liquid-solid FSI problems, two distinct methods exist for the advection of fluid interfaces: (i) solid-aware, and (ii) solid-agnostic advective schemes. Geometric VOF methods, which are inherently discontinuous interface capturing techniques, are most naturally suited to solid-aware schemes. This approach first reconstructs the updated solid interface in the domain. The gas and liquid volume fractions are then advected in the remaining domain $\Omega_f(t) = \Omega \setminus \Omega_s(t)$. By design, the gas-liquid interface does not exist inside the solid region in the geometric VOF approach; see Fig.~\ref{fig_wec_vof}.  Additionally, material triple points can be reconstructed using the geometric approach. Pathak and Raessi~\cite{pathak20163d} have combined the solid-aware geometric VOF advection scheme with the fictitious domain IB Algorithm~\ref{FDM-Chorin-algo} to solve three phase FSI problems. The solid-agnostic advective scheme, on the other hand, allows the (continuous) gas-liquid interface to pass freely through solids; see Fig.~\ref{fig_wec_ls}. The LS method is most naturally suited to this approach, and it greatly simplifies implementation. Gas-liquid interfaces within the solid region blurs the distinctions between real and fictitious fluid masses, however. For certain classes of problems, where the fluid volume is much larger than the immersed object (for example, wave energy conversion or ship hydrodynamic applications), the amount of fluid actually penetrating the immersed solid has minimal impact on the overall FSI dynamics. Several works involving large reservoir three-phase FSI problems have employed the solid-agnostic LS advection scheme~\cite{sanders2011new,bergmann2022numerical,BhallaBP2019,Bergmann2015,Khedkar2020,khedkar2022model}. If the fluid volume is significantly smaller or comparable to that of the solid, additional constraints must be imposed in the advection scheme to conserve the total volume of fluid outside the body, if the gas-liquid interface exists within the immersed object. Such constraints have recently been discussed in Khedkar et al.~\cite{khedkar2023conserving}.  
}

\section{\REVIEW{Miscellaneous}} \label{sec_miscellaneous}
\REVIEW{The purpose of this section is to discuss miscellaneous issues as well as improvements of the immersed boundary method that have been proposed in the literature. While this is not an exhaustive list (the IB literature is extremely vast), it covers some key developments.}

\subsection{High density ratio flows}

Fluid-structure interaction applications involving high density ratios ($\ge 100$) between different flowing/fluid phases are usually susceptible to numerical instabilities at high Reynolds number~\cite{nangia2019dlm,nangia2019robust,patel2017novel,pathak20163d,raessi2012consistent,desjardins2010methods}. We emphasize that numerical instabilities occur due to high density contrasts between flowing phases; adding a solid phase to the domain does not cause numerical instabilities on its own, if extended domain methods are used. For example, consider a two-phase flow in which a solid is fully immersed in a \emph{uniform/single} fluid phase. In this case the FSI solution obtained using the fictitious domain method retains stability for models involving extremely low, nearly equal, equal, and high solid-fluid density ratios~\cite{apte2009numerical,kolahdouz2021sharp}. For non-fictitious domain methods, the solution does not extend inside the structure region. Typically, such methods suffer from the ``added mass" effects of lighter structures embedded in dense fluids.

\begin{algorithm}[h!]
  \caption{HIGH DENSITY RATIO STABILIZED FICTITIOUS DOMAIN METHOD} \label{HDR-FDM-algo}
  \begin{algorithmic}
    \State \textbf{1.} Set $\rho_f^n$ using $\phi^n$. \Comment{scalar $\phi$ implicitly tracks various phases} \\

    \State \textbf{2.}  Integrate the mass conservation equation \\ \\
    $\quad \quad \D{\rho_f}{t} + \V{\nabla} \cdot (\rho_f \u_f) = 0$, \Comment{e.g. use SSP-RK3 integrator} \\ 
    to obtain $\rho_f^{n+1}$. Store the mass flux $\widetilde{\rho_f\u_f}$ used in this transport equation. \\
    
    \State \textbf{3.}  Advect the level set or VOF variable(s)  \\ 
    $\quad \quad \D{\phi}{t} + \u_f \cdot \V{\nabla} \phi  = 0$. \Comment{$ \phi^n \longrightarrow \phi^{n+1}$} \\
    
     \State \textbf{4.} Set the viscosity field based on latest $\phi^{n+1}$ or advect viscosity directly: $\; \mu_f^n \longrightarrow \mu_f^{n+1}$. \\

     \State  \textbf{5.} Solve the \emph{conservative} form of fluid momentum equation in the entire domain $\Omega$ using updated $\rho_f^{n+1}$ and $\mu_f^{n+1}$ \\ \\
    $\quad \quad  { {\rho_f^{n+1}\widehat{\u}_f - \rho_f^n\u_f^n}\over \Delta t} + \V{\nabla} \cdot (\widetilde{\rho_f\u_f} \widehat{\u}_f)= {\bm\nabla}\cdot{\widehat{\bm \sigma}} + \rho_f^{n+1}\V{g}$, \\ \\
    $\quad \quad \bm\nabla \cdot \widehat{\u}_f = \V{0}$.
    \Comment{use the same mass flux from step 2 in the momentum convective operator} \\
     
    \State \textbf{6.} Solve rigid velocity component of the solid \\ \\
    $\quad \quad \int \limits_{\mathcal{B}} {\bm \lambda}^{n+1} dV = {\bm 0} \implies {\bm U}_s^{n+1}$, 
    \Comment{external forces are assumed to be absent} \\
    
    $\quad \quad	\int \limits_{\mathcal{B}} ({\bm r} \times {\bm \lambda}^{n+1}) dV = {\bm 0} \implies {\bm \omega}_s^{n+1}$, 
    \Comment{external torques are assumed to be absent} \\
    
    $\quad \quad {\bm \lambda}^{n+1} = \rho_f^{n+1} { {{\bm U}_s^{n+1} + ({\bm \omega}_s^{n+1} \times {\bm r}) + {\bm u}_{s,\textrm{def}}^{n+1} - \widehat{\bm u}_f }\over \Delta t}$. \Comment{$\V{\lambda}$ defined on $\mathcal{B}(t)$}.  \\
    
    \State \textbf{7.} Impose velocity constraint on $\mathcal{B}(t)$ \\  \\
    $\quad \quad \rho_f^{n+1} { {{\bm u}^{n+1} - \widehat{\bm u}_f }\over \Delta t} = {\bm \lambda}^{n+1}.$  
     \Comment{Chorin-type projection}
    
  \end{algorithmic}
  \end{algorithm}

A number of approaches have been proposed to stabilize high-density ratio flows in recent years, including volume of fluid, level set, and other interface tracking methods~\cite{vaudor2017consistent,le2013monotonicity,owkes2017mass,jemison2014compressible,duret2018pressure}. A central idea of these stabilizing techniques is to discretely match the mass flux used in the momentum convective operator with the mass flux used in advecting the density field. Due to the strong coupling between two convective operators, conservative momentum and mass conservation equations are required, rather than non-conservative equations used in Equations~\ref{eq:multi-eqns1} and~\ref{eq:multi-eqns4}. We solve a single momentum and mass conservation equation by matching the fictitious fluid density to the solid density in the solid region in the multiphase fictitious domain algorithm. Using a single density field $\rho_f(\x,t)$ to represent all phases in $\Omega$, Algorithm~\ref{HDR-FDM-algo} describes the stabilized FDM approach for high density ratio flows, as implemented in Nangia et al.~\cite{nangia2019dlm}. Self-propelling bodies are considered here.

Step 1 of the algorithm involves \emph{synchronizing} the density field with the level set or VOF scalar variable $\phi$ at time level $n$. In general, several scalar variables can be used in the simulation to track different interfaces like liquid-gas or solid-liquid interfaces separately. In that case, $\phi$ in Algorithm~\ref{HDR-FDM-algo} represents the set of such variables. Synchronization of the density field with the interface tracking variable $\phi$ prevents the interface from getting diffused over time by direct advection of $\rho_f$. Step 2 of the algorithm advects the density field to the new time level $n+1$ using a discrete mass flux $\widetilde{\rho_f \u_f}$. During step 3, the scalar variable is advected to the new time level $n+1$, which in turn is used in step 4 to update the viscosity field $\mu_f^{n+1}$. Direct advection of the viscosity field is also possible, but rarely practiced. Using the updated density and viscosity fields, step 5 solves the momentum and continuity equations in conservative form. The mass and momentum advection are strongly coupled in this step by using the same discrete mass flux from step 2. For high density ratio flows, this step is crucial to maintaining numerical stability.  The rigidity constraints in the solid domain in steps 6 and 7 remain essentially the same as in Algorithm~\ref{FDM-Chorin-algo}.  
 
\subsection{Hydrodynamic force and torque evaluation}

When using the velocity forcing method, the net hydrodynamic force and torque acting on the rigid body are required to update the translational and rotational velocities. In contrast, the fictitious domain Algorithm~\ref{FDM-Chorin-algo} does not need to explicitly evaluate hydrodynamic force and torque to evolve rigid body dynamics. Hydrodynamic stress and force on the surface are still relevant, even as postprocessed information. Aerodynamic performance metrics like drag and lift coefficients can be evaluated this way. 

The fictitious domain or velocity forcing methods are typically implemented as diffuse-interface methods, where the Lagrange multipliers $\V{\Lambda}_s$ or $\V{\Lambda}_b$ enforcing rigidity constraints are smeared using regularized integral kernels. In the conventional numerical realization of the IB method, regularized kernels are used to smooth stress jumps at the interface, implying that stresses do not converge pointwise. However, as discussed in Goza et al.~\cite{goza2016accurate}, the zeroth (net hydrodynamic force) and first moment (net hydrodynamic torque) of stresses still converge at a first-order rate. It results from solving the discrete integral equation of the first-kind explicitly enforcing the Dirichlet boundary condition of Equation~\ref{eq:vel-constraint}. Even so, it is desirable to obtain smooth and pointwise convergent stresses on the body surface with immersed techniques. In the following sections, we investigate some algorithms that have been used to obtain smooth and possibly convergent hydrodynamic forces and torques at diffuse interfaces. 

\subsubsection{Smooth net hydrodynamic force and torque} 

If the objective is to obtain net hydrodynamic force and torque acting on the surface of the body, then this can be calculated easily via Lagrange multipliers. Therefore, evaluating pressure or velocity derivatives over the complex surface of the immersed body to obtain these quantities is avoidable in this situation. This was discussed in the context of the velocity forcing Algorithm~\ref{VFM-algo}, and the corresponding equations are summarized below

\begin{align}
\label{eq:net_F}
\int \limits_{\partial\Omega_s} {\bm \sigma}_f^{+} \cdot \V{\hat{n}} dS &= {D \over Dt} \left(\; \int \limits_{\Omega_s} \rho_f {\bm u}_f dV \right)  - \int \limits_{\Omega_s} {\bm \lambda} dV, \\ 
\label{eq:net_T}
\int \limits_{\partial\Omega_s} \V{r} \times ({\bm \sigma}_f^{+} \cdot \V{\hat{n}}) dS &= {D \over Dt} \left(\; \int \limits_{\Omega_s}{\bm r} \times \rho_f {\bm u}_f dV \right) - \int \limits_{\Omega_s} ({\bm r} \times {\bm \lambda}) dV.
\end{align}

Goza et al.~\cite{goza2016accurate} showed that the zeroth moment of $\V{\lambda}$ remains smooth irrespective of the width of the regularized kernel employed to smear the constraint force. Furthermore, the smoothness of the integral is not affected by kernel differentiability. For example, the authors in~\cite{goza2016accurate} demonstrated the smoothness of $\int {\bm \lambda} dV$ using  $\mathcal{C}^0$, $\mathcal{C}^1$, $\mathcal{C}^2$, $\mathcal{C}^3$, and $\mathcal{C}^\infty$ kernels --- all produced similar results. Similar conclusions can be drawn for the first moment of  the constraint force, i.e. $ \int \V{r} \times {\bm \lambda} dV$ is also smooth. The time derivative terms on the right hand side of Equations~\ref{eq:net_F} and~\ref{eq:net_T} represent the linear and angular acceleration of the body, respectively. For bodies moving with prescribed kinematics, the acceleration terms are also smooth (assuming the prescribed motion is smooth). However, when bodies are freely moving, time derivative terms generally contain spurious oscillations. These oscillations can however be easily filtered out using a moving point average operator (or any other filtering technique) without discarding physically relevant information; see Bhalla et al.~\cite{bhalla2013unified} for example cases where this technique has been applied. Therefore, computing net hydrodynamic force and torque using Equations~\ref{eq:net_F} and~\ref{eq:net_T} leads to smooth evaluation of the hydrodynamic quantities.  

The right hand side of Equations~\ref{eq:net_F} and~\ref{eq:net_T} is convenient to evaluate using a Lagrangian representation of the immersed body. For a pure Eulerian implementation of the IBM (e.g., the Brinkman penalization method), Nangia et al.~\cite{nangia2017moving} showed that by using the Leibniz integral rule, the integrals over the $\Omega_s$ region appearing in Equations~\ref{eq:net_F} and~\ref{eq:net_T}  can be shifted over a moving control volume surrounding the immersed body. The moving domain can conveniently be selected as a rectangular region that conforms to the Cartesian grid lines. The moving control volume can have a velocity different from the body it tracks. Results shown in~\cite{nangia2017moving} confirm that the net hydrodynamic force and torque calculated via the (Eulerian) moving control volume approach and via the (Lagrangian) Lagrange multiplier approach are equivalent. Moreover, the transformation preserves the smoothness of hydrodynamic quantities.

 \subsubsection{Smooth pointwise hydrodynamic force and torque}  
 
Sometimes it is also desirable to evaluate the pointwise value of hydrodynamic stresses in order to analyze the distribution of forces and torques over the surface of the body. To obtain a smooth representation of these quantities, Verma et al.~\cite{verma2017computing} recommend evaluating them over a ``lifted" surface towards the fluid side ($\partial\Omega_s^+$). If the actual surface of the body is represented by the zero-contour of a level set field, then the lifted surface can be thought of as a level set contour having a value proportional to the half-width of the regularized kernel. The results of Verma et al. demonstrate that evaluating pressure or  velocity derivatives on the lifted surface leads to smooth stress values. 

\subsubsection{Smooth pointwise constraint forces}  

In the context of Lagrange multipliers, Goza et al.~\cite{goza2016accurate} proposed an inexpensive filtering technique using redistribution of constraint forces to smooth the $\V{\Lambda}_s$ field. Furthermore, the proposed technique can be applied as a postprocessing step. By design, the technique does not affect the convergence of the original velocity or pressure field. Their results demonstrate that although the original $\V{\Lambda}_s$ field does not converge under grid refinement (due to regularized kernels), its filtered counterpart converges. We remark that $\V{\Lambda}_s$ represents jump in $\V{\sigma}_f^+$ and $\V{\sigma}_f^-$ as shown in Equation~\ref{eq:fluid-jump}, and should not be interpreted as the pointwise value of hydrodynamic force acting on the surface of the body, which was discussed in the previous section. Nevertheless, the technique proposed by Goza et al. produces smooth and convergent values of constraint forces, if these are desired --- for example to compute pointwise force field $\V{F}_s$ as described in Algorithm~\ref{FDM-algo}. 

\subsection{Sharp implementation of surface constraint force $\V{\Lambda}_s$}

Immersed body techniques are conventionally implemented numerically by smearing the surface constraint force $\V{\Lambda}_s$ using a regularized integral kernel with a non-zero compact support. As an alternative, it is possible to achieve a sharp representation of the constraint force on the fluid grid. Here we briefly discuss Kolahdouz et al.'s approach~\cite{kolahdouz2020discrete,kolahdouz2021sharp,kolahdouz2023sharp} to implement surface constraint forces in a sharp manner.

Consider a version of the velocity forcing Algorithm~\ref{VFM-step1-algo} where the constraint forces are applied only to body surface $\partial \Omega_s(t)$, such that $\V{\lambda} = \V{\Lambda}_s \delta_s$. This form of constraint force is taken  in~\cite{kolahdouz2020discrete,kolahdouz2021sharp,kolahdouz2023sharp}. The central idea behind imposing constraint force sharply is to explicitly resolve the jump condition given in Equation~\ref{eq:fluid-jump} on the fluid grid, which we re-write below taking a specific form of $\V{\Lambda}_s$ 
\begin{align}
\label{eq:lambda_s_jump}
\llbracket \V{\sigma}_f \rrbracket \cdot\hat{\bm n} &= \biggl[{\bm\sigma}_f^+ - {\bm \sigma}_f^-\biggl]\cdot\hat{\bm n} = -{\bm \Lambda}_s, \\
\label{eq:lambda_s}
{\bm \Lambda}_s(\x,t) & \approx \widetilde{\V{\Lambda}_s}(\x,t) = k(\x_s - \x) + c \left(\u_s(\x,t) - \u_f(\x,t) \right),  &  \x~\textrm{on}~\partial \Omega_s(t) 
\end{align}
in which $k$ is a spring-like penalty parameter, $c$ is a damper-like penalty parameter, and $\x_s$ and $\u_s$ are the prescribed position and velocity of the immersed surface $\partial\Omega_s(t)$. In the limit of $k \rightarrow \infty$ and/or $c \rightarrow \infty$, the constraint force $\widetilde{\V{\Lambda}_s}$ approaches the true Lagrange multiplier $\V{\Lambda}_s$. The advantage of imposing a weak constraint given by Equation~\ref{eq:lambda_s} is that it avoids a fully-implicit FDM implementation. In what follows, the discussion ignores the specific form of the constraint force. Both implicit and explicit treatments of $\V{\Lambda}_s$ are possible.

\begin{figure}[]
  \centering
   \includegraphics[scale = 0.52]{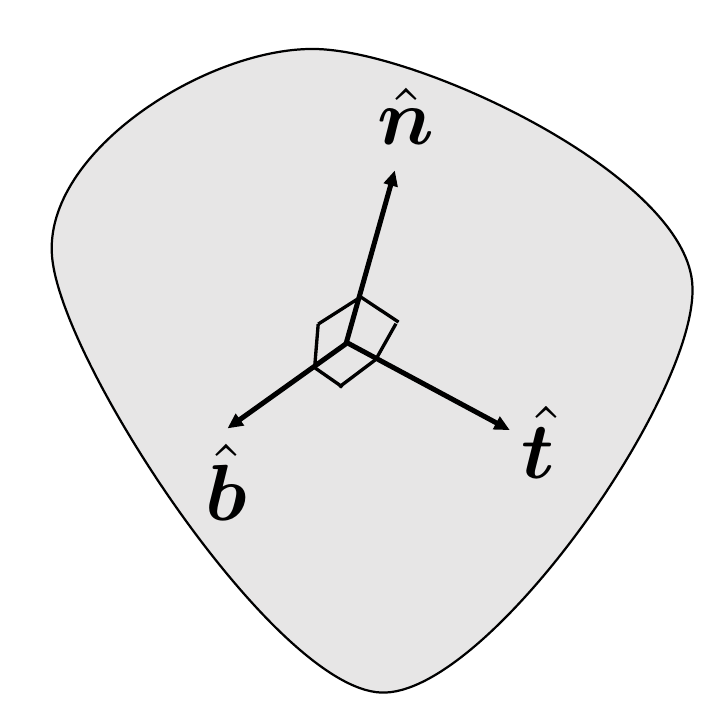}
   \caption{Local unit normal ($\hat{\V n}$) and tangent ($\hat{\V t}$ and $\hat{\V b}$) vectors of the interface. A local orthonormal system is chosen to derive the jump conditions.}
     \label{fig_normal_tangent}
\end{figure}

To proceed, we take the form of the fictitious fluid stress tensor ${\bm \sigma}_f^-$ same as that of the outside fluid stress tensor ${\bm \sigma}_f^+$. Moreover, the fluid viscosity across the interface is also taken to be the same, so that $\llbracket \mu_f \rrbracket = 0$. Since the flow is viscous, the velocity components $\left(\u_f \equiv (u,v,w)\right)$, as well as their tangential derivatives are continuous across the interface:
\begin{align}
&\llbracket u \rrbracket = \llbracket v \rrbracket = \llbracket w \rrbracket = 0, \\
\label{eq:dudt}
&\llbracket \nabla u \cdot \hat{\V t} \rrbracket = \llbracket \nabla v \cdot \hat{\V t} \rrbracket  = \llbracket \nabla w \cdot \hat{\V t} \rrbracket  = 0, \\
\label{eq:dudb}
&\llbracket \nabla u \cdot \hat{\V b} \rrbracket = \llbracket \nabla v \cdot \hat{\V b} \rrbracket  = \llbracket \nabla w \cdot \hat{\V b} \rrbracket  = 0.
\end{align}
In the above, $\hat{\V t}$ and $\hat{\V b}$ are local tangent vectors of the interface; see Fig.~\ref{fig_normal_tangent}. The incompressibility condition across both sides of the interface implies  $\llbracket \V{\nabla} \cdot \u \rrbracket  = 0$. In 
component form this can be expressed as
\begin{equation}
\label{eq:dudn}
\llbracket (\nabla u \cdot \hat{\V n},  \nabla v \cdot \hat{\V n}, \nabla w \cdot \hat{\V n}) \cdot \hat{\V n} \rrbracket  + \underbrace{\llbracket  (\nabla u \cdot \hat{\V t},  \nabla v \cdot \hat{\V t}, \nabla w \cdot \hat{\V t} ) \cdot \hat{\V t}  \rrbracket}_{= \; 0} + \underbrace{\llbracket  (\nabla u \cdot \hat{\V b},  \nabla v \cdot \hat{\V b}, \nabla w \cdot \hat{\V b}) \cdot \hat{\V b} \rrbracket}_{= \; 0} = 0.
\end{equation}
Using Equations~\ref{eq:dudt} and~\ref{eq:dudb}, the last two jump terms in the above equation vanish. In this case, Equation~\ref{eq:dudn} expresses continuity of the normal derivative of the normal component of velocity  
\begin{equation}
\label{eq:jumpdudn}
\left\llbracket  \frac{\partial \left(\u_f \cdot \hat{\V n} \right)}{\partial  \hat{\vec{n}} } \right\rrbracket = 0.
\end{equation}  

By taking the dot product of Equation~\ref{eq:lambda_s_jump} with $\hat{\V n}$, incorporating Equation~\ref{eq:jumpdudn} and the assumption $\llbracket \mu_f \rrbracket = 0$, it is straightforward to show that the discontinuity in pressure across the interface 
is given by
 \begin{align}
\label{eq:p_jump} 
\llbracket p_f(\vec{x},t) \rrbracket &=   \V{\Lambda}_s (\vec{x},t) \cdot\hat{\bm n}.  
\end{align}
Similarly, by taking the dot product of Equation~\ref{eq:lambda_s_jump} with $\hat{\V t}$ and $\hat{\V b}$, the discontinuity in the viscous stress is given as
  \begin{align}
\label{eq:viscous_stress_jump} 
\mu_f \left\llbracket \frac{\partial \left(\u_f \cdot \hat{\V t} \right)}{\partial  \hat{\vec{n}} } \right\rrbracket  =
-\V{\Lambda}_s \cdot  \hat{\V t} \quad~\textrm{and}~\quad \mu_f \left\llbracket \frac{\partial \left(\u_f \cdot \hat{\V b} \right)}{\partial  \hat{\vec{n}} } \right\rrbracket  = -\V{\Lambda}_s \cdot  \hat{\V b},~\textrm{respectively}. 
\end{align}
In component form, Equations~\ref{eq:jumpdudn} and~\ref{eq:viscous_stress_jump}  can be expressed as

\begin{equation}
\begin{pmatrix}
\hat{\V n} \\
\hat{\V t} \\
\hat{\V b}
\end{pmatrix}
\begin{pmatrix}
\llbracket \nabla u \rrbracket \\
\llbracket \nabla v \rrbracket\\
\llbracket \nabla w \rrbracket
\end{pmatrix}
\begin{pmatrix}
\hat{\V n} \\
\hat{\V t} \\
\hat{\V b}
\end{pmatrix}^\intercal
= \frac{1}{\mu_f}
\begin{pmatrix}
0 & 0 & 0 \\
-\V{\Lambda}_s \cdot  \hat{\V t} & 0 & 0 \\
-\V{\Lambda}_s \cdot  \hat{\V b} & 0 & 0
\end{pmatrix},
\end{equation}
or more explicitly as

\begin{equation}
\label{eqn_comp_visc_stress_jump}
\begin{pmatrix}
\left\llbracket  \frac{\partial u}{\partial x} \right\rrbracket & \left\llbracket  \frac{\partial u}{\partial y} \right\rrbracket & \left\llbracket  \frac{\partial u}{\partial z} \right\rrbracket \\ 
\left\llbracket  \frac{\partial v}{\partial x} \right\rrbracket & \left\llbracket  \frac{\partial v}{\partial y} \right\rrbracket & \left\llbracket  \frac{\partial v}{\partial z} \right\rrbracket \\ 
\left\llbracket  \frac{\partial w}{\partial x} \right\rrbracket & \left\llbracket  \frac{\partial w}{\partial y} \right\rrbracket & \left\llbracket  \frac{\partial w}{\partial z} \right\rrbracket 
\end{pmatrix}
= \frac{1}{\mu_f}
\begin{pmatrix}
\hat{\V n} \\
\hat{\V t} \\
\hat{\V b}
\end{pmatrix}^\intercal
\begin{pmatrix}
0 & 0 & 0 \\
-\V{\Lambda}_s \cdot  \hat{\V t} & 0 & 0 \\
-\V{\Lambda}_s \cdot  \hat{\V b} & 0 & 0
\end{pmatrix}
\begin{pmatrix}
\hat{\V n} \\
\hat{\V t} \\
\hat{\V b}
\end{pmatrix}.
\end{equation}

The jump conditions expressed in Equations~\ref{eq:p_jump} and~\ref{eqn_comp_visc_stress_jump} can be incorporated directly into the finite difference stencils used in the momentum equation, instead of regularizing the constraint force through smooth kernels. On the fluid grid, this allows for a sharp representation of surface constraint forces. We refer readers to Kolahdouz et al.~\cite{kolahdouz2020discrete,kolahdouz2021sharp,kolahdouz2023sharp} for implementation details, where the authors followed the modern immersed interface method (IIM) approach~\cite{le2006immersed,xu2006systematic} and employed generalized Taylor series expansions to incorporate the physical jump conditions into the finite difference stencils, without having to modify linear solvers. This is in contrast to the original IIM introduced by LeVeque and Li~\cite{leveque1994immersed} for elliptic PDEs with discontinuous coefficients and singular forces, which modifies the Poisson operator near the interface, and consequently requires using a different linear solver (instead of a standard Poisson solver). The original IIM technique has subsequently been extended to the incompressible Stokes~\cite{leveque1997immersed,li2007augmented} and Navier-Stokes~\cite{li2001immersed,lee2003immersed,lee2003immersed} equations, and has also been combined with level set methods~\cite{xu2006level,hou1997hybrid,sethian2000structural}. 

We remark that the aforementioned approach of implementing constraint forces sharply should not be confused with other sharp interface methods described in the literature, such as the embedded boundary formulation~\cite{yang2006embedded, kim2006immersed}, the cut-cell method~\cite{schneiders2016efficient, pogorelov2018adaptive,muralidharan2016high}, the curvilinear immersed boundary method~\cite{gilmanov2005hybrid,borazjani2008curvilinear}, and the ghost-cell immersed boundary method~\cite{tseng2003ghost}; unlike the Lagrange multiplier approach considered in this section, these methods solve fluid equations outside the immersed object and zero-out the solver solution within the fictitious fluid domain. These sharp methods do not rely on constraint forces, but impose velocity matching conditions directly, either through velocity reconstruction~\cite{gilmanov2005hybrid} or through cut-cell approaches~\cite{udaykumar2001sharp}.  

\subsection{Incorporating Neumann and Robin boundary conditions on immersed surfaces}

The fictitious domain approach to modeling fluid-structure interactions can also be applied to systems describing heat and mass transfer or chemical reactions occurring over immersed surfaces. In these systems, all three type of boundary conditions are relevant (Dirichlet, Neumann, and Robin). In contrast, only the Dirichlet boundary condition (for velocity) is required to solve the fluid-structure interaction problem, as discussed in Section~\ref{sec_unified_formulation}. This section discusses how to impose different types of boundary conditions on the surface of a solid. The transport equation for a general scalar variable $q(\x,t)$ is considered. We consider the Robin boundary condition on solid surfaces for generality. The equation system reads as 
\begin{align}
\label{eq:adv_diff}
\D{q(\x,t)}{t} +  \u_f (\x,t) \cdot\grad q(\x,t) &= \div \left[\mathcal{D}\, \grad q(\x,t) \right] +  f(\x,t), & \x~\textrm{in}~\Omega \\ 
\label{eq:adv_diff_robin_bc}
a q + b \D{q}{n} &= g(\x,t), &   \x~\textrm{on}~\partial \Omega_s(t)
\end{align}     
in which $\mathcal{D}$ is the diffusion coefficient which could be spatially or temporally varying, $f(\x,t)$ is a source/sink term, and $\u_f(\x,t)$ is the advection velocity which in general is obtained by solving the coupled (to the transport equation) fluid-structure interaction system. It is assumed that $q$ satisfies appropriate boundary conditions on the computational domain boundary $\partial \Omega$ and we omit those boundary conditions while discussing the above system of equations. 

The transport equation along with the Robin boundary condition on the solid surface can be cast into a single equation using the volume penalization (VP) approach. The volume penalized form of  transport equation reads as
\begin{equation}
\D{q}{t} + \left(1 - \chi \right) \left(\u_f\cdot\grad q\right) =   \div \left[\mathcal{D} \, \grad q(\x,t) \right] +\left(1 - \chi\right) f + \frac{\chi}{\eta}\left(a q + b \D{q}{n} - g \right), \label{eq:adv_diff_VP}
\end{equation}
in which $\eta$ is the penalization parameter, and  $\chi(\x,t)$ is the characteristic function of the body such that $\chi(\x,t) = 1$ if $\x \in \Omega_s(t) \cup \partial \Omega_s$ and $\chi(\x,t) = 0$ if $\x \notin \Omega_s(t)$. By construction, in the fluid domain where $\chi(\x,t) = 0$, the VP Equation~\ref{eq:adv_diff_VP} reverts to its original form given by Equation~\ref{eq:adv_diff}. In the limit of penalization parameter $\eta \rightarrow 0$, it can be shown that the solution to the penalized Equation~\ref{eq:adv_diff_VP} converges to the solution of the original system~\cite{bensiali2015penalization}.

The volume penalized Equation~\ref{eq:adv_diff_VP} can be solved in different ways depending upon the temporal treatment of the penalized Robin boundary condition term (the last term) of Equation~\ref{eq:adv_diff_VP}.  Brown-Dymkoski et al.~\cite{brown2014characteristic} and Hardy et al.~\cite{hardy2019penalization} treated this term explicitly to model the energy transport equation in the context of compressible flows. In contrast, Bensiali et al.~\cite{bensiali2015penalization} treated the volume penalized Robin boundary condition term implicitly and incorporated it into system of equations.  When Dirichlet boundary condition on the solid surface is desired, i.e. when $a = 1$ and $ b = 0$ in Equation~\ref{eq:adv_diff_robin_bc}, the VP transport Equation~\ref{eq:adv_diff_VP} reverts to the original Brinkman penalization formulation introduced by Angot et al.~\cite{angot1999penalization}. 

A different formulation for imposing homogeneous and spatially-varying inhomogeneous Neumann boundary conditions on $\partial \Omega_s$ ($a = 0$ and $ b = \mathcal{D}$ in Equation~\ref{eq:adv_diff_robin_bc}) is proposed by Thirumalaisamy et al.~\cite{thirumalaisamy2020,thirumalaisamy2022handling}. In their approach the authors modify the diffusion operator to impose Neumann boundary conditions. The VP transport equation for the Neumann boundary condition reads as~\cite{thirumalaisamy2020,thirumalaisamy2022handling} 
\begin{equation}
\label{eq:neumann_bc}
\D{q}{t} + \left(1 - \chi \right) \left(\u_f\cdot\grad q\right) = \div \left[\left\{\mathcal{D} \left(1 - \chi\right) + \eta \chi \right\} \grad q \right]+\left(1 - \chi\right) f + \div \left(\chi \V{\beta}\right) - \chi \div \V{\beta},
\end{equation}
in which $\V{\beta}(\x,t)$ is an arbitrary flux forcing function such that $\hat{\V{n}} \cdot \V{\beta} = g(\x,t)$ on the interface $\partial \Omega_s$. 

\begin{figure}[]
  \centering
   \includegraphics[scale = 0.35]{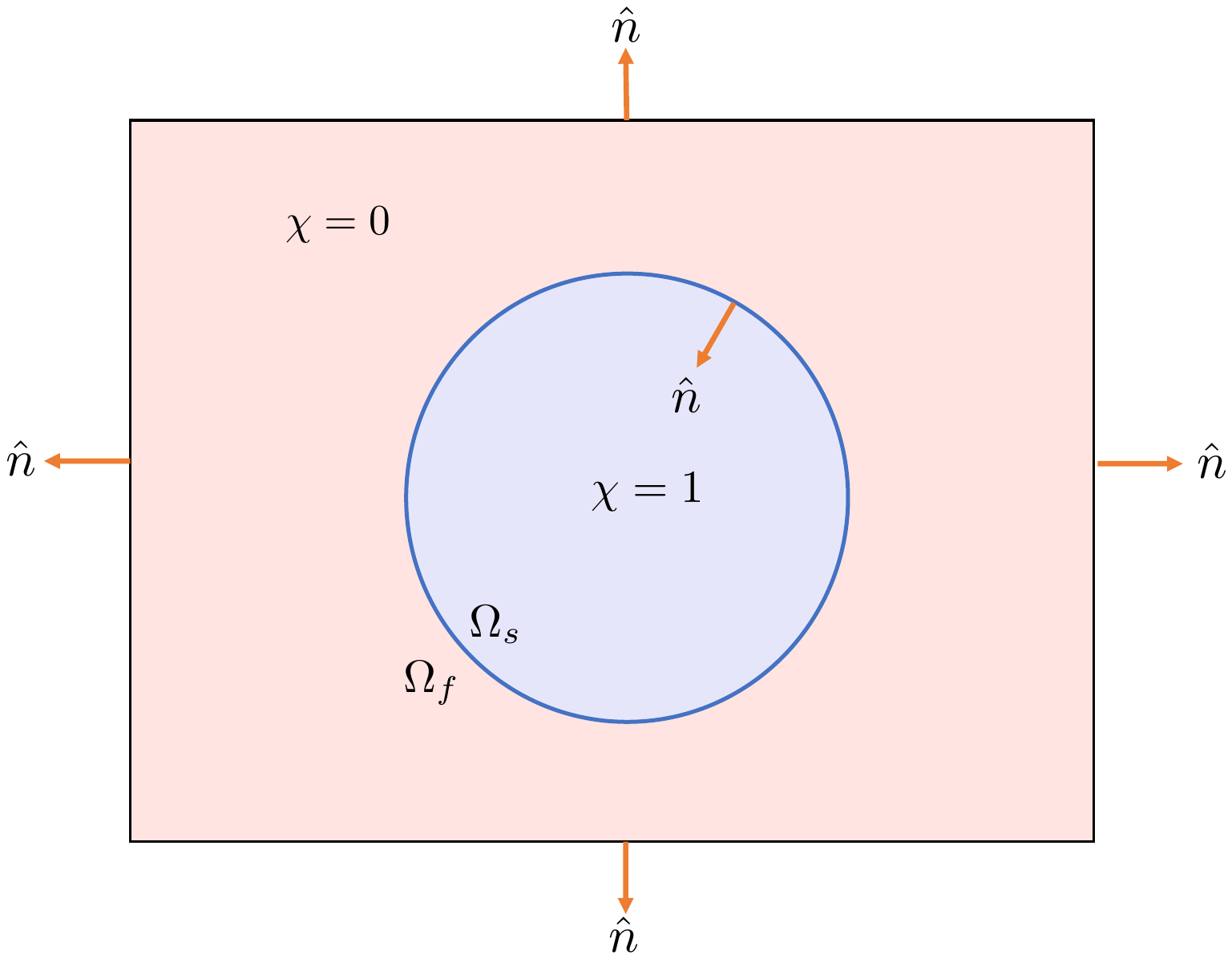}
   \caption{Imposing spatially-varying Neumann and Robin boundary conditions on an embedded interface using penalized Poisson Equation~\ref{eq:penalized_poisson}. Here $\hat{\V{n}}$ is the unit outward normal of the fluid region which is shown by filled pink region. }
     \label{fig_normal_tangent}
\end{figure}

A formal derivation of Equation~\ref{eq:neumann_bc} appears in~\cite{thirumalaisamy2022handling}. Here we present a less formal and intuitive explanation of how Equation~\eqref{eq:neumann_bc} imposes flux/Neumann boundary conditions on the immersed interface. Consider the time-independent and penalized version of Poisson equation
\begin{equation}
\label{eq:penalized_poisson}
-\div \left[\left\{\mathcal{D} \left(1 - \chi\right) + \eta \chi \right\} \grad q + \chi \V{\beta} \right] = \left(1 - \chi\right) f - \chi \div \V{\beta},
\end{equation}
which is obtained by dropping the time-derivative and convective terms from Equation~\ref{eq:neumann_bc}. Integrating Equation~\ref{eq:penalized_poisson} in the fluid region, and applying the Gauss-divergence theorem yields
\begin{align}
\label{eq:penalized_poisson_derivation}
&-\int_{\Omega_f} \left( \div \left[\left\{\mathcal{D} \left(1 - \chi\right) + \eta \chi \right\} \grad q + \chi \V{\beta} \right] \right) \; dV = \int_{\Omega_f} \left( \left(1 - \chi\right) f - \chi \div \V{\beta} \right) \; dV, \\
&-\int_{\partial \Omega} \mathcal{D} \grad q \cdot \hat{\V{n}} \; dS  -\int_{\partial \Omega_s}  \left[\eta  \grad q + \V{\beta} \right] \cdot \hat{\V{n}} \; dS =  \int_{\Omega_f} f \; dV,  \quad  [\textrm{using}~\chi = 0~\textrm{in}~\Omega_f,\textrm{ and}~\chi = 1~\textrm{in}~\Omega_s\textrm{ and on }\partial \Omega_s ]\\
&-\int_{\partial \Omega} \mathcal{D} \grad q \cdot \hat{\V{n}} \; dS  -\int_{\partial \Omega_s}  \V{\beta} \cdot \hat{\V{n}} \; dS =  \int_{\Omega_f} f \; dV. \quad \quad [\textrm{with}~\eta \rightarrow 0]
\end{align}
Therefore, in the limit of penalization parameter $\eta \rightarrow 0$, Equation~\ref{eq:penalized_poisson} converges to the solution of non-penalized Poisson equation with flux boundary condition imposed on the embedded interface.  \REVIEW{Thirumalaisamy et al. have also extended their flux boundary condition approach to impose spatially-varying Robin boundary conditions ($a = \zeta$ and $ b = \mathcal{D}$ in Equation~\ref{eq:adv_diff_robin_bc}) on  $\partial \Omega_s$. The penalized form of the Poisson equation with Robin boundary conditions on the immersed interface reads as~\cite{thirumalaisamy2022handling}
\begin{equation}
\label{eqn_vp_robin}
\zeta [\div (\chi \hat{\V{n}})  - \chi \div \hat{\V{n}}]  q  -\grad \cdot [\left\{ \mathcal{D} \left(1 - \chi \right) + \eta \chi \right \} \; \grad q] = (1- \chi) f + \div (\chi \V{\beta}) - \chi \div \V{\beta}.
\end{equation} 
Here, again the vector-valued forcing function $\V{\beta}$ satisfies the condition $\hat{\V{n}} \cdot \V{\beta} = g(\x,t)$ on the interface. A formal derivation of Equation~\ref{eqn_vp_robin} is provided in~\cite{thirumalaisamy2022handling}. In general, constructing $\V{\beta}$ is non-trivial for geometrically complex interfaces. The authors in~\cite{thirumalaisamy2022handling} present a signed distance function-based numerical technique for constructing $\V{\beta}$ for complex interfaces.} 


\subsection{``Leakage"} \label{sec_leakage}

With diffuse-interface immersed boundary methods, constraint forces are smeared using regularized integral kernels. A typical isotropic kernel function couples the structure to the fluid degrees of freedom on both sides of the interface. The use of isotropic kernels leads to spurious feedback forces and internal flows typically observed in diffuse-interface IB models. Such types of spurious flows within a solid are referred to as ``leaks"  in the IB literature. The most common recommendation to make a structure ``water-tight" using immersed body techniques is to impose $\u_s = \u_f$ on $\mathcal{B}(t)$, instead of imposing it on $\partial \Omega_s(t)$ alone. While volumetric imposition of a no-slip constraint reduces spurious flows within a solid, they are not entirely eliminated. Moreover, the ``best" relative grid spacing between the Lagrangian and Eulerian mesh that minimizes leakage needs to be found empirically. For example, for Peskin's IBM Algorithm~\ref{IBM-algo}, the recommended ratio of Lagrangian grid spacing and Eulerian cell size is between half and one-third~\cite{griffith2012volume} (i.e., two or three Lagrangian points in each coordinate direction inside an Eulerian grid cell), whereas for the fictitious domain Algorithm~\ref{FDM-Chorin-algo} it is one~\cite{curet2010versatile,bhalla2013unified}. Furthermore, if the volumetric version of the no-slip constraint is imposed using the fully-implicit fictitious domain approach, then the size of the problem (i.e., the number of degrees of freedom related to constraint force calculation) and the computational cost of the solution is increased significantly, especially when three-dimensional bodies are considered~\cite{kallemov2016immersed,balboa2017hydrodynamics}. 

Recently, Bale et al.~\cite{bale2020} proposed a \emph{one-sided} direct forcing IBM using kernel functions constructed via \emph{moving least squares} (MLS)~\cite{wendland2001local,mirzaei2012generalized,liu1997moving}. In their approach, the kernel functions effectively couple structural degrees of freedom to fluid variables on only one side of the fluid-structure interface. Their results demonstrate that one-sided IB kernels reduce spurious feedback forcing and internal flows typically observed in diffuse-interface IB models for relatively high Reynolds number flows. We briefly outline the MLS approach to constructing one-sided kernel functions next.

In the moving least squares method, a quasi-interpolant of the type

\begin{equation}
\cP f(\x) = \sum_i^N f(\x_i) \psi_i(\X),      
\end{equation}
is sought.  In the above equation  $\cP$ is the interpolation operator, $\f = [f(\x_1), \ldots, f(\x_N)]^T$ are given data at $N$ \emph{interpolation} points, and $\X$ is the \emph{evaluation} point. 
The moving least squares method computes \emph{generating functions} $\psi_i(\X) = \{\psi(\x_i,\X)\}$ subject 
to the polynomial reproduction constraints
\begin{equation}
\InterpSumi p(\x_i) \psi_i(\X) = p(\X), \quad \text{for all } p \in \Pi_{d}^s,  \label{eqn_poly_reproduce}
\end{equation}
in which $\Pi_{d}^s$ is the space of s-variate polynomials of total degree at most $d$.
The size of polynomial basis set is $m = \frac{(s+d) !}{d!}$. The polynomial reproduction constraints 
correspond to \emph{discrete moment conditions} for the function $\psi_i(\X)$. In the matrix form, Equation~\ref{eqn_poly_reproduce} is written as
\begin{equation}
\cA \V{\Psi} (\X) = \V P(\X),  \label{eqn_linsys}
\end{equation}
in which the entries of the polynomial matrix $\cA \in \mathbb{R}^{m \times n}$ are the values of the basis functions at the data point 
locations, $\cA_{ij} = p_i(\x_j), i = 1, \ldots,m, j = 1, \ldots, N$, and the right-hand side vector $\V P = [p_1, \ldots, p_m]^T$
contains the values of the polynomials at the evaluation point $\X$. The unknown generating function 
vector $\V{\Psi} = [\psi_1, \ldots, \psi_N]^T$ is obtained by solving a \emph{least squares problem}.
Generally $N \gg m$, and therefore, the underdetermined system of equations~\ref{eqn_linsys}  is solved in a weighted least squares sense with Lagrange multipliers $\V{\Upsilon} (\X)$ to enforce the reproducing conditions. 
With $W(\x_i,\X)$ denoting a positive weight function that decreases with distance from the evaluation point $\X$,
the Lagrange multipliers are obtained by inverting the Gram matrix $\cG(\X)  = \cA \cW \cA^T \in \mathbb{R}^{m \times m}$  

\begin{equation}
\cG (\X) \V{\Upsilon} (\X) = \V P(\X).  \label{eqn_gram}
\end{equation}
Here $\cW(\X) =  \text{diag}\left( \V W \right)$ is the diagonal matrix containing weights, and 
$\V W = \left(W(\x_1,\X), \dots, W(\x_N,\X) \right)$ is the main diagonal of $\cW(\X)$. 
The symmetric positive-definite Gram matrix (assuming polynomial basis is linearly independent)
has weighted  $L^2$ inner products of the polynomials as its entries
\begin{equation}
\cG_{jk}(\X) = \left< p_j, p_k \right>_{W(\X)} = \InterpSumi p_j(\x_i) p_k(\x_i) W(\x_i, \X), \quad j,k = 1,\ldots,m. \label{eqn_gram_mat}
\end{equation}  
The generating functions are obtained using the Lagrange multipliers $\V{\Upsilon} (\X)$ 

\begin{equation}
 \psi(\x_i, \X) = \psi_{i}(\X)  = W(\x_i,\X) \InterpSumj \Upsilon_j(\X) p_j(\x_i), \quad i = 1, \ldots, N. \label{eqn_psi}
\end{equation}
In matrix form, Equation~\ref{eqn_psi} can be written as
\begin{equation}
 \V{\Psi} (\X)  =  \cW (\X) \odot \vcL(\X), \label{eqn_matrix_psi}
\end{equation}
in which $\vcL(\X) = \cA^T \V{\lambda}(\X) $ and $\odot$ indicates component-wise product of two matrices.

\begin{figure}
	\centering
\includegraphics[scale=0.4]{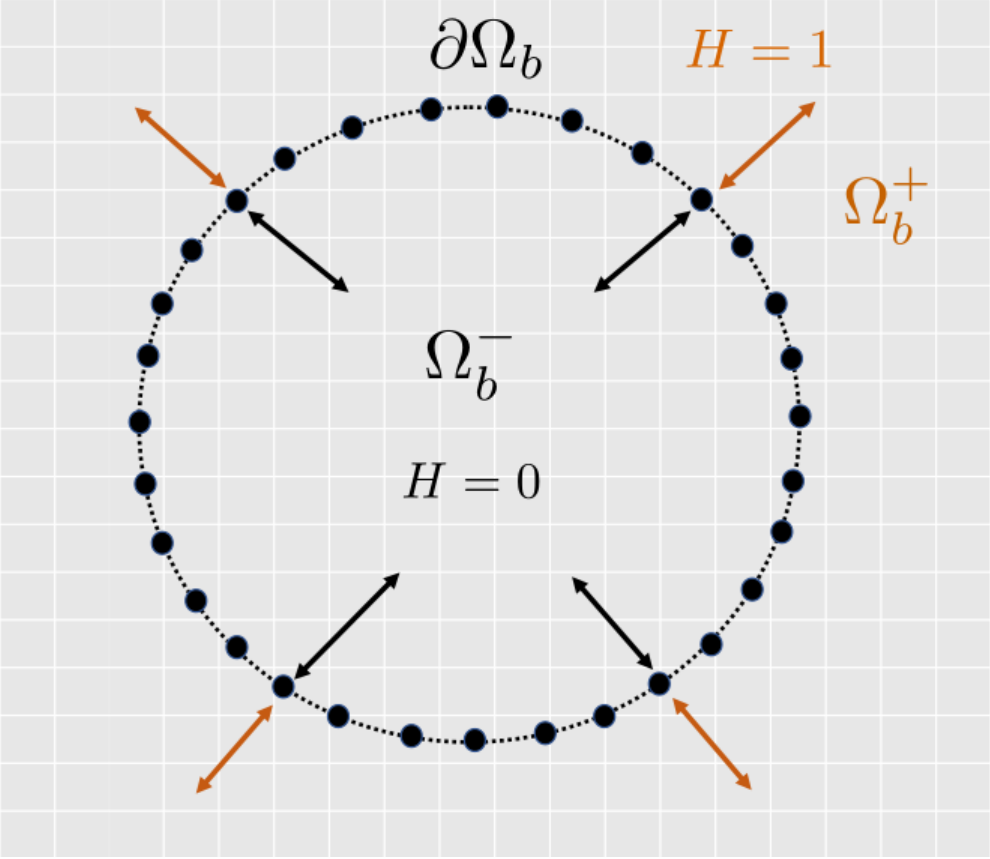}	 
\caption{The Eulerian computational domain $\Omega$ is partitioned into an inner region $\Omega_b^{-}$, and an outer region $\Omega_b^{+}$  based on the location of the interface $\partial \Omega_b$. Lagrangian marker points $[\bullet]$ on the interface interpolate velocity from and spread force to  $\Omega_b^{-}$ and $\Omega_b^{+}$ regions separately. The Heaviside function $H(\x)$ is defined to be $H = 1$ in $\Omega_b^{+}$, and $H = 0$ in $\Omega_b^{-}$. In this example,  $\Omega_b^{+}$ is the interaction region for the Lagrangian marker points.}
	\label{fig:MLSSchematic}
\end{figure}

Now, to obtain the one-sided IB kernel $\psi^{\circ}(\x_i,\X)$ using the MLS procedure, define the Heaviside function $H(\x)$
\begin{align}
H(\x) &= 
\begin{cases} 
       0,  & \Omega_b^{-},\\
       1  &  \Omega_b^{+} \cup \partial \Omega_b,  \label{eq_heaviside_body}
\end{cases}
\end{align}
to select the appropriate side of the interface $\partial \Omega_b$ from where the Eulerian velocity field is interpolated, and the Lagrangian force is spread; the domain of influence for the evaluation/Lagrangian point $\X$ is the spatial region having value of $H(\x) = 1$ (see Fig.~\ref{fig:MLSSchematic}). 
Next,  multiply the standard/unrestricted weights by the Heaviside function to obtain the restricted weights $W_{\text{MLS}}$
\begin{equation}
W_{\text{MLS}}(\x_i,\X) = W(\x_i,\X)H(\x_i) \quad \text{s.t.} \quad  W_{\text{MLS}}(\x_i,\X) = 0 \quad  \forall \quad \x_i \in \Omega_b^{-}.  \label{eqn_wmls}
\end{equation}
It can be easily verified that by using $W_{\text{MLS}}$ in Equation~\ref{eqn_psi}, $\psi^{\circ}(\x_i,\X) = 0\; \forall\;  \x_i \in \Omega_b^{-}$.  

\REVIEW{
\subsection{Dense particle-laden flows}

\begin{figure}
	\centering
\includegraphics[scale=0.4]{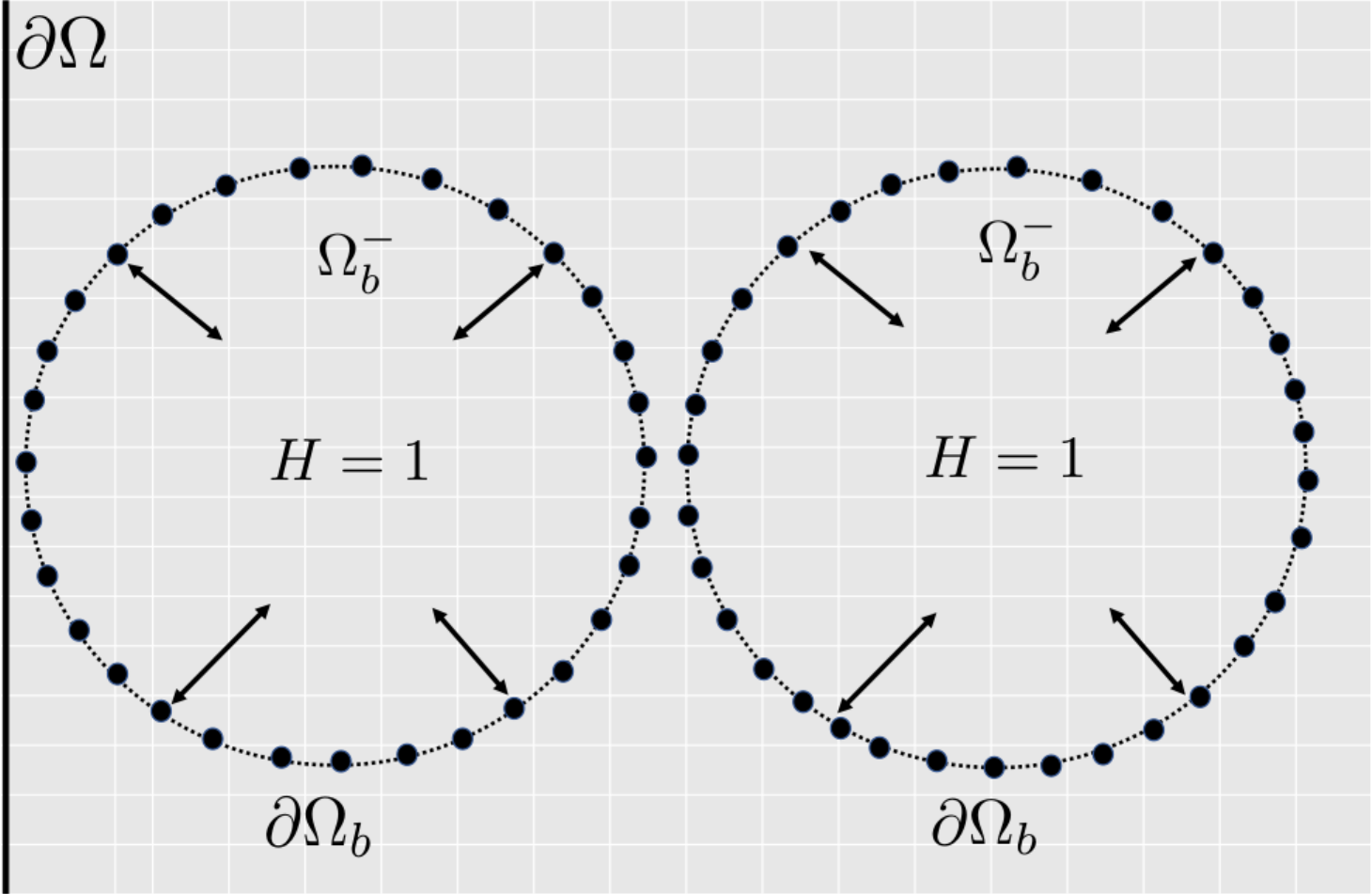}	 
\caption{One-sided IB kernel approach for dense particle laden flows. Lagrangian marker points $[\bullet]$ on the interface interpolate velocity from and spread force to  particle interior region $\Omega_b^{-}$. The Heaviside function $H(\x)$ is defined to be $H = 1$ in $\Omega_b^{-}$, and $H = 0$ in $\Omega_b^{+}$. Here, $\Omega_b^{-}$ is the interaction region for the Lagrangian marker points.}
	\label{fig:MLSSchematic2}
\end{figure}

Flows laden with particles are ubiquitous in natural and industrial processes, such as dust storms, river sediment transport, and the dispersion of particles in chemical reactors. In order to understand particle behavior in flows, numerical simulations are critical.  It is essential to consider the presence of other particles and their interaction with the fluid phase when simulating particle-laden flows using CFD. Particle-laden flows can be modeled using the IB method, which allows direct numerical simulation (DNS). The IB kernel support can go inside nearby particles' interiors or outside the computational domain when two particles are close together, or reach near the computational boundary. The accuracy of IB simulations is significantly reduced in both scenarios. In order to avoid particle collisions or close proximity, artificial repulsion forces are typically used in IB simulations~\cite{glowinski1999distributed,Patankar2000}. 

Recently, Ch\'{e}ron  et al.~\cite{cheron2023hybrid,cheron2023drag} presented a promising approach in which they employed the one-sided IB kernels of Bale et al.~\cite{bale2020} (that were also discussed in Section~\ref{sec_leakage}) to simulate dense particle-laden flows. In their approach, the particle interior is used as the Eulerian domain for interpolating velocity and spreading Lagrangian constraint forces. Fig.~\ref{fig:MLSSchematic2} shows a schematic illustration of the approach. Thus, one-sided IB kernels are not affected by particle gaps. Results in~\cite{cheron2023hybrid,cheron2023drag} demonstrate that a one-sided IB approach achieves better results than classical IBM frameworks, especially when particles are close together or near domain walls.      

}


\section{Future challenges} \label{sec_future}
In the diffuse-interface IB method the constraint forces are smeared on the background grid using regularized integral kernels, whereas in the sharp interface IB method the jump conditions (across the interface) arising from the constraint forces are explicitly resolved on the fluid grid. The former approach is easier to implement in practice as it requires only kernel function weights to transfer quantities between Eulerian and Lagrangian grids. In contrast, the latter approach requires computational geometry constructs to identify the Eulerian grid points where jump conditions need to be applied. This leads to a substantial amount of bookkeeping for complex geometries and in three spatial dimensions. An alternate approach could be to generate kernel function weights that automatically satisfy the desired jump conditions. Such an approach would combine the ease of implementation of diffuse-interface methods to the improved spatial accuracy offered by the sharp interface IB methods.  One possible way to achieve this objective would be to extend the moving least squares approach~\cite{wendland2001local,mirzaei2012generalized,liu1997moving} to incorporate jump conditions along with the usual polynomial reproduction constraints to produce such kernel function weights. 

Making fully implicit fictitious domain methods efficient and competitive against their explicit or semi-implicit counterparts is still an active area of research. Stiff  and nonlinear material models for elastic bodies, and dense mobility matrices for rigid bodies provide significant challenges to the monolithic FSI solvers. A possible strategy to reduce the runtime cost of monolithic solvers would be to leverage modern GPUs to perform factorization of small but dense systems in batch mode. Another possibility would be develop efficient physics-based preconditioners that leverage (semi-) analytical solutions --- even if for some selective FSI cases. 

Immersed boundary methods have traditionally been used in the context of single phase flows. Their use for multiphase flows is relatively recent. Along the line of multiphase flow modeling, a challenging yet exciting area of application for IB methods is additive manufacturing, in which a solid body undergoes phase transition. Diffuse-interface immersed boundary methods would be more appropriate for these applications, as the transitioning phase boundary is physically ``mushy". Recently this approach has been used to model melting and solidification problems in the presence of a passive gas phase~\cite{thirumalaisamy2023low}. Extending this methodology to model the simultaneous occurrences of melting, solidification, evaporation and condensation, which happen during various metal manufacturing processes (e.g., welding, selective laser melting/sintering, thin film deposition), would be a worthwhile research avenue.  

Immersed methods coupled with AI approaches have the potential to lead to efficient computational tools. Finally, the extended domain strong form immersed formulation might be explored in turbulence and particulate closure models.





\section*{Acknowledgements} NAP acknowledges support from NSF grants OAC 1450374 and 1931372. APSB acknowledges support from NSF awards OAC 1931368 and CBET CAREER 2234387.



\begin{flushleft}
 \bibliography{ReviewBibliography}
\end{flushleft}

\end{document}